\documentclass[pdftex,11pt]{article}
\usepackage{url}
\usepackage{graphicx}
\usepackage[usenames]{color}
\usepackage{subfigure}
\usepackage{latexsym}
\usepackage{hyperref}
\usepackage{pslatex}
\usepackage{amsmath}
\usepackage{amsfonts}

\renewcommand{\S}{\mathbf{S}}
\newcommand{\x}{\mathbf{x}}
\newcommand{\y}{\mathbf{y}}
\newcommand{\U}{\mathbf{U}}
\newcommand{\M}{\mathbf{M}}
\newcommand{\A}{\mathbf{A}}
\newcommand{\D}{\mathbf{D}}
\newcommand{\Q}{\mathbf{Q}}
\newcommand{\R}{\mathbf{R}}
\newcommand{\T}{\mathbf{T}}
\renewcommand{\L}{\mathbf{L}}
\newcommand{\I}{\mathbf{I}}
\newcommand{\C}{\mathbf{C}}
\newcommand{\e}{\mathbf{e}}
\renewcommand{\b}{\mathbf{b}}
\renewcommand{\u}{\mathbf{u}}
\renewcommand{\v}{\mathbf{v}}
\newcommand{\p}{\mathbf{p}}
\newcommand{\q}{\mathbf{q}}
\newcommand{\rr}{\mathbf{r}}
\renewcommand{\t}{\mathbf{t}}
\newcommand{\g}{\mathbf{g}}
\renewcommand{\l}{\mathbf{l}}
\renewcommand{\ll}{\mathbf{z}}

\definecolor{DarkRed}{rgb}{0.55,0.00,0.00}
\definecolor{DeepPink2}{rgb}{0.93,0.07,0.54}

\hypersetup{
pdftitle={},
pdfauthor={Julien Langou}, 
} 

\title{
Translation and modern interpretation of Laplace's
Th\'eorie Analytique des Probabilit\'es, pages 505-512, 516-520.
}

\author{Julien Langou}

\begin{document}
\maketitle

\begin{abstract}
The text of Laplace, \textit{Sur l'application du calcul des probabilit\'es \`a
la philosophie naturelle,} (Th\'eorie Analytique des Probabilit\'es.
Troisi\`eme \'Edition.  Premier Suppl\'ement), 1820, is quoted in the context
of the Gram-Schmidt algorithm. We provide an English translation of Laplace's
manuscript (originally in French) and interpret the algorithms of Laplace in a
contemporary context. The two algorithms given by Laplace computes the mean and
the variance of two components of the solution of a linear statistical model.
The first algorithm can be interpreted as {\em reverse square-root-free
modified Gram-Schmidt by row} algorithm on the regression matrix. The second
algorithm can be interpreted as the {\em reverse square-root-free Cholesky} algorithm.
\end{abstract}

\section{Introduction}
\label{sec:introduction}

This translation work is inspired by the one of Pete Stewart~\cite{Stewart.95}
who translates from Latin to English the ``{\em Theoria Combinationis
Observationum Erroribus Minimis Obnoxiae}'' of Gauss in the SIAM book
``{\em Theory of the Combination of Observations Least Subject to Errors, Part
One, Part Two, Supplement.}'' Stewart translates 101 original pages of
Gauss, and he also provides an important contribution (28 pages) to place the work of
Gauss in a historical framework. This manuscript is a more modest contribution.
I translate thirteen pages and explain the relation of Laplace's algorithm
with our contemporary algorithms. I would like to thanks Pete Stewart to have
inspired me by his work.  I also would like to thank {\AA}ke Bj\"orck for
giving me my first version of Laplace's manuscript back in 2004 and Serge Gratton for useful comments on an early draft of the manuscript.

The goal of Laplace is to compute the mass of Jupiter (or Saturn) from a system
of normal equations provided by Bouvard and from this same system to compute
the distribution of error in the solution assuming a normal distribution of the
noise on the observations. The parameters of the noise distribution are not known. 
Laplace explains how to compute the standard deviation of two
variables of a linear statistical model. His algorithm can be interpreted as
performing the Cholesky factorization of the normal equations and then compute
the two standard deviations from the Cholesky factor. A second method used by Laplace to
justify the first is to perform a QR factorization of the regression matrix and
compute the standard deviation from the R factor. Laplace was performing the QR
factorization through the modified Gram-Schmidt algorithm.

Laplace did not know what a factorization was, nor a matrix. I interpret his
result through factorizations but certainly do not claim that Laplace invented
all this.

The first method which Laplace introduces consists in successivelly projecting the
system of equations orthogonally to a column of the observation matrix. This action
eliminates the associated variable from the updated system.
Ultimately, Laplace eliminates all the variables but the one of interest in the
linear least squares problem, which eliminates all the columns but one in the
observation matrix. Laplace is indeed introducing the main idea behind the
Gram-Schmidt algorithm (successive orthogonal projections.)
Laplace gives an example on a $s$--by--$6$ system.
Once the observation matrix is reduced to a vector column, Laplace is able to
relate the standard deviation of the variable of interest to the norm of this vector
column and the norm of the residual vector.

 While Laplace could have stopped here and performed the modified
Gram-Schmidt algorithm onto the original overdetermined system, he explains how
to compute the norm of the projected column of observations of interest directly from the
normal equations. He observes that, if he performs a Cholesky factorization of
the normal equations, the last coefficient computed will be equal to the norm
of the last column orthogonally projected successivelly to the span of the remaining columns.
In the mean time, Laplace
observes that this method (Cholesky) provides a way to get the value of
the solution from the normal equations. Laplace also generalizes this approach
to more than one variable. 

Laplace has used the modified Gram-Schmidt algorithm as a tool to derive the
Cholesky algorithm on the normal equations. Laplace did not interpret his
results with orthogonality, in particular, he did not observe that, after
orthogonal projection with respect to the last column, all the remaining
projected columns were made orthogonal to that column. The orthogonal
projections are interpreted as elimination conserving orthogonality with the
residual.  Laplace correctly explains and observes the property that all the
remaining projected columns, after elimination/projection, are orthogonal to
the residual of the least squares problem and that the residual vector is conserved.

Laplace then uses his Cholesky algorithm to solve two $6$--by--$6$ systems of
normal equations given to him by the French astronom Bouvard to recompute the mass of
Jupiter and Saturn, the originality of the work consists in assessing the
reliability of these computations by estimating the standard deviation of the
distribution of error in the solution.

In Section~\ref{sec:background}, I set up the background for Laplace's work.
This background is briefly recalled in Section 1 of Laplace's manuscript. I
chose not to translate this Section directly. It is failry hard to read indeed
and I have preferred to explain it and refer to the equations in it.  In
Section~\ref{sec:translation}, I provide a translation from French to English
of Laplace's \textit{Sur l'application du calcul des probabilit\'es
\`a la philosophie naturelle,} (Th\'eorie Analytique des Probabilit\'es.
Troisi\`eme \'Edition.  Premier Suppl\'ement), 1820. I have translatted
Sections 2 and 5 which represent pp.505-512 and pp.516-520.

This manuscript of Laplace is quoted in the book of
Farebrother~\cite[Chap.4]{Farebrother.88} and the book of
Bj\"orck~\cite[p.61]{Bjorck.96}). In both books, the authors claim
that Laplace is using the modified Gram-Schmidt algorithm.

The text is available in French from the Biblioth\`eque Nationale de France website
\footnote{See 
\texttt{http://gallica.bnf.fr/ark:/12148/bpt6k775950}.
For Section 2, type $505$ in the box \textit{Aller Page}.
For Section 5, type $516$ in the box \textit{Aller Page}.}.
There is one typo (in the pages we are translatting). On p.517, the seventh equation should read
$-668486",70 = -13208350 z + 413134432 z' - \textit{151992,0} z''' - 34876,7 z^{iv}.$

We present some terminology used by Laplace.  The overdetermined system of
equations is named: \textit{les \'equations g\'en\'erales de condition des
\'el\'ements}.  The Linear Least Squares method is named: \textit{la m\'ethode
la plus avantageuse}.  The \textit{poids} (=weight) $P$ of a normal
distribution is related to the standard deviation $\sigma$ by $P =\left( 2 \sigma^2 \right)^{-1}$.

\section{Background}
\label{sec:background}

The main goal of this manuscript is to provide a translation of Laplace's
algorithmic contribution.  However to put things into context, we start in this
section with some background and notations.

\subsection{Covariance matrix of the regression coefficients of a linear statistical model}

Laplace considers the classical linear statistical model
$$
   \A\x  = \b^*+ \e,
   ~\A \in \mathbb{R}^{s\times n},
   ~\b^* \in \mathbb{R}^s,
$$
where $\e$ is a vector of random errors.
having normal distribution and we 
will denote
$\sigma_e$ the standard deviation of $\|\e\|$.
In statistical language, the matrix $\A$ is referred to as the regression matrix
and the unknown vector $\x$ is called the vector of regression coefficients.
In our context, the matrix $\A$ is full rank.
If $\e=0$ (there is no error in the data), then we denote the solution of the consistent overdetermined linear system of equations as $\x^\star$,
$$ \A\x^\star = \b^* .$$

Given a vector of $s$ observations $\b$,
Laplace considers the linear estimate $\x$ given by the solution of the 
linear least squares method.
We call $\e'$ the residual of the linear least squares solution
$$\e' \equiv  \A\x - \b .$$
In Laplace terms, (see, e.g., [p.501] last sentence), by the
{\it conditions de la m\'ethode la plus avantageuse}, we have
$$
   \sum p^{(i)}\epsilon'^{(i)} = 0, \quad
   \sum q^{(i)}\epsilon'^{(i)} = 0, \quad \ldots;
$$
where $p^{(i)}$ is the element (i,1) of $\A$, $q^{(i)}$ is the element (i,2) of $\A$, $\ldots$,
and $\epsilon^{(i)}$ is the element i of $\e'$.
In other words,
$$
   \e' =  \A\x - \b \perp \textmd{Span}~(\A).
$$
Several other definitions for $\x$, the linear least squares solution,  are possible. We give two more equivalent definitions
  $$\x \textmd{ is such that } \left\| \A\x - \b \right\| = \min_{\y} \left\| \A\y - \b\right\|,
 \quad\textmd{or, equivalently,}\quad  \x \textmd{ is such that } \A^T\A \x = \A^T \b.$$
We define $\u=(u,u',\ldots)$ the random vector which represents the error between the vector $\x^{*}$ and the linear estimator $\x$.
From p.501 to p.504, Laplace derives the formula of the joint distribution of the random variables $u$, $u'$, $\ldots$ 

{\em We} know that the joint distribution of the multivariate centered normal variable $\u$ is  proportional to
  $$ \exp\left( -\frac{1}{2} \u^T\C^{-1}\u \right),$$
where $\C$ is the covariance matrix of $\u$.  In our case we have 
  $$ \C \equiv \sigma_b^2 \left( \A^T\A \right) ^{-1} ,$$
therefore, the joint distribution of $\u$ is  proportional to
\begin{equation}\label{eq:l}
   \exp\left( -\frac{1}{2\sigma_b^2} \u^T \A^T\A \u \right).
\end{equation}
In practice one does not know $\sigma_b^2$ and we therefore rely on an unbiased estimate, for example
 $$\frac{1}{s-n}\left\|\b -\A\x\right\|^2.$$
Now if we approximate $s-n$ with $s$, we obtain that
the random vector $\u$ follows a multivariate normal distribution with covariance matrix
  $$ \frac{ 1 }{s}\| \b - \A\x\|^2\left( \A^T\A \right) ^{-1} =
   \frac{ 1 }{s}\| \e' \|^2\left( \A^T\A \right) ^{-1} .$$
So that the joint distribution of the random variables $u$, $u'$, $\ldots$ 
is proportional to
  $$ \exp\left( - \frac{ s }{ 2 \|\e'\|^2 } \u^T \A^T \A \u \right) .$$
This result is given in term of the variable $\v \equiv \u/\sqrt{s}$ by Laplace.
Laplace states that the  joint distribution of the random variables $v$, $v'$, $\ldots$ 
is proportional to (see first formula top of p.504)
  $$ \exp\left( - \frac{ \sum \left( p^{(i)}v + q^{(i)}v' + \ldots \right)^2 } {2 \sum \epsilon'^{(i)2} } \right) .$$
Note that
  $$ \exp\left( - \frac{ s }{ 2 \|\e'\|^2 } \u^T \A^T \A \u \right)
  = \exp\left( - \frac{ \sum \left( p^{(i)}v + q^{(i)}v' + \ldots \right)^2 } {2 \sum \epsilon'^{(i)2} } \right) .$$
Therefore, Laplace's framework fits our standard linear statistical model framework.

\subsection{Laplace's algorithm to compute of the standard deviation of one variable of a linear statistical model}
\subsubsection{Laplace and the modified Gram Schmidt algorithm}

The background is now half set. We have a linear statistical model to which we
seek a regression vector $\x$ through the linear least squares method and we know
that the covariance matrix of the regression vector is given by the matrix $\C$.
Laplace wants to compute only the first variable, $z$, of the regression vector, $\x$.
He also seeks the standard deviation, $\sigma_z=\sigma_u$, of this variable. 

From p.504 to p.505, Laplace explains that the modified Gram-Schmidt process
applied to the matrix $\A$ enables him to find the standard deviation of the first variable.
Laplace applies the modified Gram-Schmidt in a backward manner, that is, he projects the columns
$1$ to $n-1$ orthogonally to the span of the column $n$ and obtains the matrix $\A_1$,
then, working from the updated matrix $\A_1$, he projects the columns $1$ to $n-2$
orthogonally to the span of the column $n-1$, etc.

The {\em reverse modified Gram-Schmidt by row algorithm} on the matrix $\A$ is formally given as follows.

Following Laplace, we name the columns of $\A$
$$ \A = \left( \p, \q, \rr, \ldots \t, \g, \l\right). $$

First step is to project column $1$ to $n-1$ of $\A$ orthogonally to the span of its column $n$, so we define
\begin{eqnarray}
\nonumber&& \p_1 = \left( \I - \frac{\l~\l^T}{\|\l\|^2} \right) \p, \\
\nonumber&& \q_1 = \left( \I - \frac{\l~\l^T}{\|\l\|^2} \right) \q, \\
\nonumber&& \vdots \\
\nonumber&& \t_1 = \left( \I - \frac{\l~\l^T}{\|\l\|^2} \right) \t, \\
\nonumber&& \g_1 = \left( \I - \frac{\l~\l^T}{\|\l\|^2} \right) \g.
\end{eqnarray}
This defines the $s$-by-$(n-1)$ matrix 
$ \A_1 = \left( \p_1, \q_1, \rr_1, \ldots \t_1, \g_1\right). $

Second step is to project column $1$ to $n-2$ of $\A_1$ orthogonally to the span of its column $n-1$, so we define
\begin{eqnarray}
\nonumber&& \p_2 = \left( \I - \frac{\g_1\g_1^T}{\|\g_1\|^2} \right) \p_1, \\
\nonumber&& \q_2 = \left( \I - \frac{\g_1\g_1^T}{\|\g_1\|^2} \right) \q_1, \\
\nonumber&& \vdots \\
\nonumber&& \t_2 = \left( \I - \frac{\g_1\g_1^T}{\|\g_1\|^2} \right) \t_1.
\end{eqnarray}
This defines the $s$-by-$(n-2)$ matrix 
$ \A_2 = \left( \p_2, \q_2, \rr_2, \ldots \t_2\right). $

At the end of step $n-2$, we have computed the $s$-by-$2$ matrix $\A_{n-2}$.
The step $n-1$ consists in projecting the first column of $\A_{n-2}$ orthogonally to the span of its second column
\begin{eqnarray}
\nonumber&& \p_{n-1} = \left( \I - \frac{\q_{n-2}\q_{n-2}^T}{\|\q_{n-2}\|^2} \right) \p_{n-2}.
\end{eqnarray}

Nowadays we are use to describe the modified Gram-Schmidt the other way around:
project orthogonally to column 1, then column 2, etc.  In either case, we note
that we need to order our variables correctly.  With Laplace's method (reverse
modified Gram-Schmidt), we will see that it is crucial to have the variable of
interest ordered first. (And ordered last in the case of forward modified
Gram-Schmidt.)

Forward modified Gram-Schmidt generates a QR factorization of the matrix $\A$,
that is, we compute
$$ 
  \A = \Q\R,
   \textmd{ where }
   \Q \textmd{ is }s\textmd{--by--}n\textmd{ with orthonormal columns and }
   \R \textmd{ is }n\textmd{--by--}n\textmd{ upper triangular.}
$$
(Without loss of generality, we will impose the diagonal elements of $\R$ to be positive.)
On the other hand, 
reverse modified Gram-Schmidt generates a QL factorization 
of the matrix $\A$, that is, we compute
$$ 
  \A = \Q\L,
   \textmd{ where }
   \Q \textmd{ is }s\textmd{--by--}n\textmd{ with orthonormal columns and }
   \L \textmd{ is }n\textmd{--by--}n\textmd{ lower triangular.}
$$
(Without loss of generality, we will impose the diagonal elements of $\L$ to be positive.)

We note that Laplace does not generate the matrix $\Q$.
Laplace generates the matrix $\T$ defined as
\begin{equation}\label{eq:T}
 \T \equiv \left( \p_{n-1}, \q_{n-2}, \rr_{n-3}, \ldots \t_{2}, \g_{1}, \l\right). 
\end{equation}
If we normalize the columns of $\T$, we will obtain $\Q$.
Laplace applies what we could call the {\em reverse square-root-free modified Gram-Schmidt by row} algorithm.
If we define
\begin{equation}\label{eq:D}
\D_M = \texttt{diag} \left(\| \p_{n-1}\|^{2}, \|\q_{n-2}\|^{2}, \|\rr_{n-3}\|^{2},
           \ldots \|\t_{2}\|^{2}, \|\g_{1}\|^{2}, \|\l\|^{2}\right)
 = \T^T \T.
\end{equation}
then we have 
$$ \T = \Q \D_M^{1/2}.$$
So that we also have the factorization
\begin{equation}\label{eq:repres}
 \A = \T \left(  \D_M^{-1/2} \L \right).
\end{equation}
The matrix $\left(  \D_M^{-1/2} \L \right)$ is lower triangular with ones on the diagonal.

\begin{center}
\begin{tabular}{ccc}
QR factorization && reverse square-root-free QR factorization \\\\
\begin{picture}(220,120)(0,0)
\scriptsize

\put(  2,  2){\line(+1,+0){ 56}}
\put( 58,  2){\line(+0,+1){106}}
\put( 58,108){\line(-1,+0){ 56}}
\put(  2,108){\line(+0,-1){106}}
\put( 20, 80){\mbox{\normalsize$\A$}}

\put( 65, 80){\mbox{\large$=$}}

\put( 82,  2){\line(+1,+0){ 56}}
\put(138,  2){\line(+0,+1){106}}
\put(138,108){\line(-1,+0){ 56}}
\put( 82,108){\line(+0,-1){106}}
\put(100, 80){\mbox{\normalsize$\Q$}}

\put(145, 80){\mbox{$\bullet$}}

\put(151,100){\mbox{$\sqrt{\bullet}$}}
\put(206, 48){\mbox{$\sqrt{\bullet}$}}
\put(208, 54){\line(-1,+1){ 40}}

\put(172,108){\line(+1,+0){ 46}}
\put(218,108){\line(+0,-1){ 46}}
\put(218, 62){\line(-1,+1){ 46}}
\put(205, 90){\mbox{\normalsize$\R$}}

\end{picture}
&~~~~~&
\begin{picture}(220,60)(0,0)
\scriptsize

\put(  2,  2){\line(+1,+0){ 56}}
\put( 58,  2){\line(+0,+1){106}}
\put( 58,108){\line(-1,+0){ 56}}
\put(  2,108){\line(+0,-1){106}}
\put( 20, 80){\mbox{\normalsize$\A$}}

\put( 65, 80){\mbox{\large$=$}}

\put( 82,  2){\line(+1,+0){ 56}}
\put(138,  2){\line(+0,+1){106}}
\put(138,108){\line(-1,+0){ 56}}
\put( 82,108){\line(+0,-1){106}}
\put(100, 80){\mbox{\normalsize$\Q\D_M^{1/2}$}}

\put(145, 80){\mbox{$\bullet$}}

\put(160, 52){\line(+1,+0){ 46}}
\put(160, 98){\line(+0,-1){ 46}}
\put(206, 52){\line(-1,+1){ 46}}
\put(156,104){\mbox{$1$}}
\put(211, 52){\mbox{$1$}}
\put(205, 60){\line(-1,+1){ 40}}
\put(162, 57){\mbox{\normalsize$\D_M^{-1/2}\L$}}

\end{picture}
\end{tabular}
\end{center}

\subsubsection{A standard relation between the standard deviation of the last variable of a statistical model and the QR factorization of the regression matrix} 
\label{subsubsec:zz}

The standard deviation of the variable $i$ of $\u$ is given by the
square-root of the entry $(i,i)$ of the covariance matrix $\C$. If we are
interested in the standard deviation $\sigma_u$ of the first variable of $\u$, $u$, we need to be able to compute the entry $(1,1)$ of
$\left( \A^T\A \right) ^{-1}$. We outline below a standard way to compute this quantity.

Once we have the QL factorization of $\A$, we write
   $$ \A^T\A = \L^T \Q^T \Q \L = \L^T \L $$
therefore
   $$
      \sigma_{u}
    = \sqrt{ \textmd{entry}(1,1)\textmd{ of }\C }
    = \sigma_b \sqrt{ \textmd{entry}(1,1)\textmd{ of }( \A^T\A )^{-1} }
    = \sigma_b \sqrt{\textmd{entry}(1,1)\textmd{ of }\L^{-1}\L^{-T}}.
   $$
And, so using the fact that the matrix $L$ is lower triangular, we have
\begin{equation}
\label{eq:mii}
      \sigma_{u} = \frac{1}{\ell_{1,1}}\sigma_b.
\end{equation}
We can prove that 
   $$ \ell_{1,1} = \| \p_{n-1}  \|, $$
(where $ \p_{n-1} $ is the vector obtained at the last step of reverse modified Gram-Schmidt algorithm), so that we obtain
that the
marginal probability density function of the first variable $z$ is
\begin{equation}
\label{eq:abcd}
   \exp\left(- \frac{ 1 } {2\sigma_b^2} \| \p_{n-1}\|^2  u ^2 \right),
\end{equation}
and if we use  the fact that
 $\frac{1}{s}\left\|\e'\right\|^2$ can be used as an approximation of 
an unbiased estimate of $\sigma_b^2$, we obtain that the
marginal probability density function of the first variable $z$ is
\begin{equation*}
     \exp\left( - \frac{1}{2} \frac{s}{\| \e' \|^2}  \| \p_{n-1}  \|^2 u^2 \right).
\end{equation*}
This formula is assessed by Laplace on top p.505.
We read:\\
``{\it This exponential becomes
$$ \exp\left(-Pu^2\right),$$
where
$$ P = \frac{s\sum p_{n-1}^{(i)2}}{2\sum \varepsilon'^{(i)2}} .$$
$u$ being the error of the random variable $z$, $P$ is what I called the \textit{poids} (weight) of this value.}''

The \textit{poids} is related to the standard deviation $\sigma$ with
$ P = \left(2\sigma^2\right)^{-1}.$
The term \textit{poids} was chosen by Laplace for the following reason (see p.499):\\

\begin{minipage}{13cm}
``{\it 
   la probabilit\'e d\'ecroit avec rapidit\'e quand il [le poids] augmente, en
   sorte que le r\'esultat obtenue p\`ese, si je puis ainsi dire, vers la
   v\'erit\'e, d'autant que ce module est plus grand.}''
\end{minipage}\\

  which gives\\

\begin{minipage}{13cm}
``{\it the probability quickly decreases with it [le poids] increases, so that
the result weights, if I can says so, towards the truth as much as this modulus
is larger.}''
\end{minipage}\\

 Other reasons are given in the same paragraph.

\subsubsection{Laplace's derivation of the standard deviation of the last variable from the QL factorization of the regression matrix}

The overall strategy of Laplace to compute $\sigma_{u}$ is well-known nowadays. How did Laplace derive it in the first place?
Starting from the  fact that the joint density function of $u$, $u'$, $\ldots$, $u^{(n)}$ is proportional to
  $$ \exp\left( - \frac{ 1 } {\sigma_b^2} \| \A u \|^2\right) ,$$
(see~Equation(\ref{eq:l}),
Laplace is interested in computed a function proportional to the
marginal probability density function of the first variable $u$, $\sigma_u$, that is Laplace wants to compute
a function of the variable $u$ proportional to
$$
    \int_{u'=-\infty}^{+\infty}
    \int_{u''=-\infty}^{+\infty} \ldots
    \int_{u^{(n)}=-\infty}^{+\infty}
  \exp\left( - \frac{ 1 } {\sigma_b^2} \| \A \u \|^2\right) du' du''\ldots du^{(n)}
$$

Laplace proposes to proceed by steps. First we will seek a function proportional to
the joint density function of $u$, $u'$, $\ldots$, $u^{(n-1)}$;
then we will seek a function proportional to
the joint density function of $u$, $u'$, $\ldots$, $u^{(n-2)}$, etc.
we will eventually end up with a function proportional to the
marginal probability density function of the first variable $u$.

To perform the first step, we therefore need to compute a
function of the variables $u$, $u'$, $\ldots$, $u^{(n-1)}$ proportional to
$$
    \int_{u^{(n)}=-\infty}^{+\infty}
  \exp\left( - \frac{ 1 } {\sigma_b^2} \| \A \u \|^2\right)  du^{(n)}.
$$

Laplace observes that, (Pythagorean theorem),
$$
\| \A\u \|^2
= 
\| \left( I - \frac{\l~\l^T}{\|\l\|^2} \right) \A\u \|^2
+
\| \frac{\l~\l^T}{\|\l\|^2}  \A\u  \|^2,
$$
and so
\begin{eqnarray}
\nonumber 
\| \A \left(\begin{array}{c}u\\u'\\\vdots\\u^{(n-1)}\\ u^{(n)}\end{array} \right) \|^2
& = &
\| \A_1 \left(\begin{array}{c}u\\u'\\\vdots\\u^{(n-1)}\end{array} \right) \|^2
+
\| \frac{\l}{\|\l\|^2} \l^T \A \left(\begin{array}{c}u\\u'\\\vdots\\u^{(n-1)}\\ u^{(n)}\end{array} \right) \|^2,\\
\nonumber
& = &
\| \A_1 \left(\begin{array}{c}u\\u'\\\vdots\\u^{(n-1)}\end{array} \right) \|^2
+
\|\l\|^2 \left(
 u^{(n)} +
 \frac{1}{\|\l\|^2}
 \l^T \left(\p,\q,\ldots,\t,\g\right) \left(\begin{array}{c}u\\u'\\\vdots\\u^{(n-1)}\end{array} \right)\right)^2
\end{eqnarray}

The joint density function of $u$, $u'$, $\ldots$, $u^{(n)}$ is therefore proportional to
$$
\exp\left(- \frac{ 1 } {\sigma_b^2} \| \A_1 \left(\begin{array}{c}u\\u'\\\vdots\\u^{(n-1)}\end{array} \right) \|^2 \right)\cdot
\exp\left(- \frac{ 1 } {\sigma_b^2}
\|\l\|^2 \left(
 u^{(n)} +
 \frac{1}{\|\l\|^2}
 \l^T \left(\p,\q,\ldots,\t,\g\right) \left(\begin{array}{c}u\\u'\\\vdots\\u^{(n-1)}\end{array} \right)
\right)^2\right).
$$
(This latter equation corresponds to the seventh equation on p.504.)

As previously explained, the first step of Laplace's derivation consists in
integrating this last term for $u^{(n)}$ ranging from $-\infty$ to $+\infty$ in
order to obtain a function proportional to the joint density function of $u$,
$u'$, $\ldots$, $u^{(n-1)}$. So let us do this. We write
\begin{eqnarray}
\nonumber
&&
\int_{u^{(n)}=-\infty}^{+\infty}
\exp\left(- \frac{ 1 } {\sigma_b^2} \| \A_1 \left(\begin{array}{c}u\\u'\\\vdots\\u^{(n-1)}\end{array} \right) \|^2 \right)\cdot
\exp\left(- \frac{ 1 } {\sigma_b^2}
\|\l\|^2 \left(
 u^{(n)} +
 \frac{1}{\|\l\|^2}
 \l^T \left(\p,\q,\ldots,\t,\g\right) \left(\begin{array}{c}u\\u'\\\vdots\\u^{(n-1)}\end{array} \right)
\right)^2\right)
du^{(n)}\\
\nonumber
&=&
\exp\left(- \frac{ 1 } {\sigma_b^2} \| \A_1 \left(\begin{array}{c}u\\u'\\\vdots\\u^{(n-1)}\end{array} \right) \|^2 \right)\cdot
\int_{u^{(n)}=-\infty}^{+\infty}
\exp\left(- \frac{ 1 } {\sigma_b^2}
\|\l\|^2 \left(
 u^{(n)} +
 \frac{1}{\|\l\|^2}
 \l^T \left(\p,\q,\ldots,\t,\g\right) \left(\begin{array}{c}u\\u'\\\vdots\\u^{(n-1)}\end{array} \right)
\right)^2\right)
du^{(n)}.
\end{eqnarray}
The second term is of the form
$$
\int_{u^{(n)}=-\infty}^{+\infty} \exp\left(- \mu \left( u^{(n)} + g( u, u',\ldots,u^{(n-1)}) \right)^2 \right) du^{(n)}.
$$
We note that this term is independent of the variables $ u, u',\ldots,u^{(n-1)}$.
Therefore we can remove this term from the previous equation and conclude that 
the joint density function of $u$,
$u'$, $\ldots$, $u^{(n-1)}$ is proportional to
\begin{equation}
\label{eq:ijk}
\exp\left(- \frac{ 1 } {\sigma_b^2} \| \A_1 \left(\begin{array}{c}u\\u'\\\vdots\\u^{(n-1)}\end{array} \right)\| ^2 \right).
\end{equation}

Continuing the process, we end up with the fact that 
the marginal probability density function of the first variable $u$
is proportional to
$$ \exp\left(- \frac{ 1 } {\sigma_b^2} \| \p_{n-1}\|^2 \cdot u ^2 \right).$$
We recover Equation~(\ref{eq:abcd}) also given on top p.505.

From this equation, Laplace deduce that, to compute the standard deviation of $u$, he needs to compute $\| \p_{n-1} \|^2$.
While it is clear that he (or Bouvard) can work on $\A$ and perform the modified Gram Schmidt algorithm, Laplace 
finds it easier to work on the normal equations. Quoting Laplace:\\

\begin{minipage}{13cm}
``{\it Mais il est plus simple d'appliquer le proc\'ed\'e dont nous venons de faire usage 
aux equations finales qui d\'eterminent les \'el\'ements, pour les r\'eduire \`a une seule, ce qui donne
une m\'ethode facile de r\'esoudre ces \'equations}.''\\
\end{minipage}\\




which means\\

\begin{minipage}{13cm}
``{\it But it is easier to apply the method we have just used to the final equations which define the variables,
in order to reduce it to a single, which gives a convenient way of solving these equations.}''
\end{minipage}\\

Therefore the next question that needs to be answered is: how can we compute  $\|
\p_{n-1} \|^2$ from $\A^T\A$ without accessing $\A$? This question is the matter
of Section 2 of Laplace's treatise from p.505 to p.512.  An example of
application of the technique is proposed in the Section 5 of the same
manuscript from p.516 to p.520.  We provide a translation of these two parts in the next section.

If we consider the QL factorization of the matrix $\A$ given by
$$ \A =
\left( \begin{array}{cc}  \Q_1 & \q \end{array} \right)
\left(
\begin{array}{cc}
\L_1 & \\
\ll^T& \alpha
\end{array}
\right),
$$
$\C_1$, the covariance matrix of joint normal distribution of the variables $u$, $u_1$,$\ldots$, $u_{n-1}$,
is 
$$ 
\C_1 = 
 \sigma_b^2
 \left( \A_1^T\A_1 \right) ^{-1} = 
 \sigma_b^2
 \left( \L_1^T\L_1 \right) ^{-1}.
$$
We can derive this relation from Laplace's analysis for example.
Another way to derive, this result is to remember
that the
covariance matrix $\C_1$ of the joint normal distribution of the variable $u$,
$u_1$,$\ldots$, $u_{n-1}$ is the $(n-1)$-by-$(n-1)$ block of the
covariance matrix $\C$ of the joint normal distribution of the variable $u$,
$u_1$,$\ldots$, $u_{n-1}, u_n$.
So if we write
\begin{eqnarray}
\nonumber
\sigma_b^2
\left( \A^T\A \right)^{-1} & = &
\sigma_b^2
\left(
\left(
\begin{array}{cc}
\L_1^T &\ll \\
& \alpha
\end{array}
\right)
\left(
\begin{array}{cc}
\L_1 & \\
\ll^T& \alpha
\end{array}
\right) \right)^{-1}\\
\nonumber
&= &
\sigma_b^2
\left(
\begin{array}{cc}
 \L_1^T\L_1 + \ll\ll^T & \ll \alpha \\
\alpha \ll^T & \alpha^2     
\end{array}
\right)^{-1}\\
\nonumber
&= & \sigma_b^2
\left(
\begin{array}{cc}
 \left( \L_1^T \L_1 \right)^{-1} & -\frac{1}{\alpha}\left(\L_1^T\L_1\right)^{-1}\ll \\
 -\frac{1}{\alpha}\ll^T\left(\L_1^T\L_1\right)^{-1} & \frac{1}{\alpha^2}\left( \ll^T\left(\L_1^T\L_1\right)^{-1}\ll+1\right) 
\end{array}
\right),
\end{eqnarray}
we see that the the $(n-1)$-by-$(n-1)$ block is $\sigma_b^2 \left(
\L_1^T\L_1 \right) ^{-1}=\C_1$.
This is two ways to explain a standard result.

\section{Translation}
\label{sec:translation}
We now present a translation of Laplace's text.
We proceed by couple of pages. First page gives the French version. Second page gives the translatted version.
We recall that the notation $S$ stands for $\sum$.
%
%
\newpage
\setcounter{equation}{0}
\renewcommand\theequation{\arabic{equation}}
\indent
2. Reprenons l'\'equation g\'en\'erale de condition, et, pour plus de simplicit\'e, bornons-la aux six \'el\'ements
$z$, $z'$, $z''$, $z'''$, $z^{iv}$, $z^v$; elle devient alors
\begin{equation}\label{eq:orig:laplace1}
\varepsilon^{(i)}
= p^{(i)} z
+ q^{(i)} z'
+ r^{(i)} z''
+ t^{(i)} z'''
+ \gamma^{(i)} z^{iv}
+ \lambda^{(i)} z^{v}
- \alpha^{(i)}.
\end{equation}
En la multipliant par $\lambda^{(i)}$ et r\'eunissant tous les produits semblables, on aura
\begin{displaymath}
 S \lambda^{(i)} \varepsilon^{(i)} = z S \lambda^{(i)} p^{(i)} +  z' S \lambda^{(i)} q^{(i)} + \ldots - S \lambda^{(i)} \alpha^{(i)} ,
\end{displaymath}
le signe int\'egral $S$ s'\'etendant \`a toutes les valeurs de $i$, depuis $i=0$ jusqu'\`a $i=s-1$, $s$ \'etant le nombre 
des observations employ\'ees. Par les conditions de la m\'ethode la plus avantageuse, on a $S \lambda^{(i)} \varepsilon^{(i)} = 0$;
l'\'equation pr\'ec\'edente donnera donc
\begin{displaymath}
 z^{v} =
 - z^{iv} \frac{ S \lambda^{(i)} \gamma^{(i)}}{ S \lambda^{(i)2} }
 - z'''   \frac{ S \lambda^{(i)} t^{(i)}}{ S \lambda^{(i)2} }
 - z''    \frac{ S \lambda^{(i)} r^{(i)}}{ S \lambda^{(i)2} }
 - z'     \frac{ S \lambda^{(i)} q^{(i)}}{ S \lambda^{(i)2} }
 - z      \frac{ S \lambda^{(i)} p^{(i)}}{ S \lambda^{(i)2} }
 + \frac{ S \lambda^{(i)} \alpha^{(i)}}{ S \lambda^{(i)2} }.
\end{displaymath}
Si l'on substitue cette valeur dans l'\'equation~(\ref{eq:orig:laplace1}) et si l'on fait
\begin{eqnarray}
 \nonumber  \gamma_1^{(i)} & = & \gamma^{(i)} - \lambda^{(i)} \frac { S \lambda^{(i)} \gamma^{(i)} }{ S \lambda^{(i)2} },\\
 \nonumber  t_1^{(i)} & = & t^{(i)} - \lambda^{(i)} \frac { S \lambda^{(i)} t^{(i)} }{ S \lambda^{(i)2} },\\
 \nonumber  r_1^{(i)} & = & r^{(i)} - \lambda^{(i)} \frac { S \lambda^{(i)} r^{(i)} }{ S \lambda^{(i)2} },\\
 \nonumber  q_1^{(i)} & = & q^{(i)} - \lambda^{(i)} \frac { S \lambda^{(i)} q^{(i)} }{ S \lambda^{(i)2} },\\
 \nonumber  p_1^{(i)} & = & p^{(i)} - \lambda^{(i)} \frac { S \lambda^{(i)} p^{(i)} }{ S \lambda^{(i)2} },\\
 \nonumber  \alpha_1^{(i)} & = & \alpha^{(i)} - \lambda^{(i)} \frac { S \lambda^{(i)} \alpha^{(i)} }{ S \lambda^{(i)2} },
\end{eqnarray}
on aura
\begin{equation}\label{eq:orig:laplace2}
\varepsilon^{(i)}
= p_1^{(i)} z
+ q_1^{(i)} z'
+ r_1^{(i)} z''
+ t_1^{(i)} z'''
+ \gamma_1^{(i)} z^{iv}
- \alpha_1^{(i)};
\end{equation}
par ce moyen, l'\'el\'ement $z^v$ a disparu des \'equations de condition que repr\'esente l'\'equation~(\ref{eq:orig:laplace2}).
En multipliant cette \'equation par $\gamma_1^{(i)}$ et r\'eunissant tous les produits semblables, en observant ensuite que l'on a 
\begin{displaymath}
S \gamma_1^{(i)} \varepsilon^{(i)} = 0
\end{displaymath}
en vertu des \'equations 
\begin{displaymath}
0 = S \lambda^{(i)} \varepsilon^{(i)}, \quad
0 = S \gamma^{(i)} \varepsilon^{(i)}
\end{displaymath}
que donnent les conditions de la m\'ethode la plus avantageuse, on aura
\begin{displaymath}
0
= z      S \gamma_1^{(i)} p_1^{(i)} 
+ z'     S \gamma_1^{(i)} q_1^{(i)} 
+ z''    S \gamma_1^{(i)} r_1^{(i)} 
+ z'''   S \gamma_1^{(i)} t_1^{(i)} 
+ z^{iv} S \gamma_1^{(i)2}
- S \gamma_1^{(i)} \alpha_1^{(i)};
\end{displaymath}
d'o\`u l'on tire
\begin{displaymath}
 z^{iv} =
 - z'''   \frac{ S \gamma_1^{(i)} t_1^{(i)}}{ S \gamma_1^{(i)2} }
 - z''    \frac{ S \gamma_1^{(i)} r_1^{(i)}}{ S \gamma_1^{(i)2} }
 - z'     \frac{ S \gamma_1^{(i)} q_1^{(i)}}{ S \gamma_1^{(i)2} }
 - z      \frac{ S \gamma_1^{(i)} p_1^{(i)}}{ S \gamma_1^{(i)2} }
 + \frac{ S \gamma_1^{(i)} \alpha_1^{(i)}}{ S \gamma_1^{(i)2} }.
\end{displaymath}
%
%
\newpage
\setcounter{equation}{0}
\renewcommand\theequation{\arabic{equation}}
\indent
2. We consider again the overdetermined system of equations, and, for the sake of simplicity,
we restrict it to the six elements $z$, $z'$, $z''$, $z'''$, $z^{iv}$, $z^v$; it
then becomes
\begin{equation}\label{eq:trad:laplace1}
\varepsilon^{(i)}
= p^{(i)} z
+ q^{(i)} z'
+ r^{(i)} z''
+ t^{(i)} z'''
+ \gamma^{(i)} z^{iv}
+ \lambda^{(i)} z^{v}
- \alpha^{(i)}.
\end{equation}
Multiplying by $\lambda^{(i)}$ and grouping all similar products, we have
\begin{displaymath}
 S \lambda^{(i)} \varepsilon^{(i)} = z S \lambda^{(i)} p^{(i)} +  z' S \lambda^{(i)} q^{(i)} + \ldots - S \lambda^{(i)} \alpha^{(i)} ,
\end{displaymath}
the integral sign $S$ ranging for all the values of $i$, from $i=0$ to $i=s-1$, $s$ being the number of observations.
By the conditions of \textit{la m\'ethode la plus avantageuse}, we have $S \lambda^{(i)} \varepsilon^{(i)} = 0$;
the former equation consequently gives
\begin{displaymath}
 z^{v} =
 - z^{iv} \frac{ S \lambda^{(i)} \gamma^{(i)}}{ S \lambda^{(i)2} }
 - z'''   \frac{ S \lambda^{(i)} t^{(i)}}{ S \lambda^{(i)2} }
 - z''    \frac{ S \lambda^{(i)} r^{(i)}}{ S \lambda^{(i)2} }
 - z'     \frac{ S \lambda^{(i)} q^{(i)}}{ S \lambda^{(i)2} }
 - z      \frac{ S \lambda^{(i)} p^{(i)}}{ S \lambda^{(i)2} }
 + \frac{ S \lambda^{(i)} \alpha^{(i)}}{ S \lambda^{(i)2} }.
\end{displaymath}
If we replace this value in Equation~(\ref{eq:trad:laplace1}) and if we perform
\begin{eqnarray}
 \nonumber  \gamma_1^{(i)} & = & \gamma^{(i)} - \lambda^{(i)} \frac { S \lambda^{(i)} \gamma^{(i)} }{ S \lambda^{(i)2} },\\
 \nonumber  t_1^{(i)} & = & t^{(i)} - \lambda^{(i)} \frac { S \lambda^{(i)} t^{(i)} }{ S \lambda^{(i)2} },\\
 \nonumber  r_1^{(i)} & = & r^{(i)} - \lambda^{(i)} \frac { S \lambda^{(i)} r^{(i)} }{ S \lambda^{(i)2} },\\
 \nonumber  q_1^{(i)} & = & q^{(i)} - \lambda^{(i)} \frac { S \lambda^{(i)} q^{(i)} }{ S \lambda^{(i)2} },\\
 \nonumber  p_1^{(i)} & = & p^{(i)} - \lambda^{(i)} \frac { S \lambda^{(i)} p^{(i)} }{ S \lambda^{(i)2} },\\
 \nonumber  \alpha_1^{(i)} & = & \alpha^{(i)} - \lambda^{(i)} \frac { S \lambda^{(i)} \alpha^{(i)} }{ S \lambda^{(i)2} },
\end{eqnarray}
we have
\begin{equation}\label{eq:trad:laplace2}
\varepsilon^{(i)}
= p_1^{(i)} z
+ q_1^{(i)} z'
+ r_1^{(i)} z''
+ t_1^{(i)} z'''
+ \gamma_1^{(i)} z^{iv}
- \alpha_1^{(i)};
\end{equation}
by this technique, the element $z^v$ has disappeared from the system of equations represented by Equation~(\ref{eq:trad:laplace2}).
Multiplying this equation by $\gamma_1^{(i)}$, grouping all similar products, and observing that we have
\begin{displaymath}
S \gamma_1^{(i)} \varepsilon^{(i)} = 0
\end{displaymath}
from the equations
\begin{displaymath}
0 = S \lambda^{(i)} \varepsilon^{(i)}, \quad
0 = S \gamma^{(i)} \varepsilon^{(i)}
\end{displaymath}
given by the conditions of \textit{la m\'ethode la plus avantageuse}, we have
\begin{displaymath}
0
= z      S \gamma_1^{(i)} p_1^{(i)} 
+ z'     S \gamma_1^{(i)} q_1^{(i)} 
+ z''    S \gamma_1^{(i)} r_1^{(i)} 
+ z'''   S \gamma_1^{(i)} t_1^{(i)} 
+ z^{iv} S \gamma_1^{(i)2}
- S \gamma_1^{(i)} \alpha_1^{(i)};
\end{displaymath}
from which we draw
\begin{displaymath}
 z^{iv} =
 - z'''   \frac{ S \gamma_1^{(i)} t_1^{(i)}}{ S \gamma_1^{(i)2} }
 - z''    \frac{ S \gamma_1^{(i)} r_1^{(i)}}{ S \gamma_1^{(i)2} }
 - z'     \frac{ S \gamma_1^{(i)} q_1^{(i)}}{ S \gamma_1^{(i)2} }
 - z      \frac{ S \gamma_1^{(i)} p_1^{(i)}}{ S \gamma_1^{(i)2} }
 + \frac{ S \gamma_1^{(i)} \alpha_1^{(i)}}{ S \gamma_1^{(i)2} }.
\end{displaymath}

%
%
\newpage
\setcounter{equation}{2}
\renewcommand\theequation{\arabic{equation}}
\noindent
Si l'on substitue cette valeur dans l'\'equation~(\ref{eq:orig:laplace2}) et si l'on fait
\begin{eqnarray}
 \nonumber  t_2^{(i)} & = & t_1^{(i)} - \gamma_1^{(i)} \frac { S \gamma_1^{(i)} t_1^{(i)} }{ S \gamma_1^{(i)2} },\\
 \nonumber  r_2^{(i)} & = & r_1^{(i)} - \gamma_1^{(i)} \frac { S \gamma_1^{(i)} r_1^{(i)} }{ S \gamma_1^{(i)2} },\\
 \nonumber  q_2^{(i)} & = & q_1^{(i)} - \gamma_1^{(i)} \frac { S \gamma_1^{(i)} q_1^{(i)} }{ S \gamma_1^{(i)2} },\\
 \nonumber  p_2^{(i)} & = & p_1^{(i)} - \gamma_1^{(i)} \frac { S \gamma_1^{(i)} p_1^{(i)} }{ S \gamma_1^{(i)2} },\\
 \nonumber  \alpha_2^{(i)} & = & \alpha_1^{(i)} - \gamma_1^{(i)} \frac { S \gamma_1^{(i)} \alpha_1^{(i)} }{ S \gamma_1^{(i)2} },
\end{eqnarray}
on aura
\begin{equation}\label{eq:orig:laplace3}
\varepsilon^{(i)}
= p_2^{(i)} z
+ q_2^{(i)} z'
+ r_2^{(i)} z''
+ t_2^{(i)} z'''
- \alpha_2^{(i)}.
\end{equation}
En continuant ainsi, on parviendra \`a une \'equation de la forme 
\begin{equation}\label{eq:orig:laplace4}
\varepsilon^{(i)} = p_5^{(i)} z - \alpha_5^{(i)}.
\end{equation}
Il r\'esulte du n$^{\textmd{o}}$20 du Livre II que, si la valeur de $z$ est d\'etermin\'ee par cette \'equation et que $u$ soit l'erreur de cette valeur, 
la probabilit\'e de cette erreur est
\begin{displaymath}
\sqrt{ \frac{ s S p_5^{(i)2} }{ 2 S \varepsilon'^{(i)2} \pi } } 
e^{- \frac{ s S p_5^{(i)2} }{ 2 S \varepsilon'^{(i)2} } u^2 }, 
\end{displaymath}
$ S \varepsilon'^{(i)2} $ \'etant la somme des carr\'es des restes des \'equations de condition, 
lorsqu'on y a substitu\'e les \'el\'ements d\'etermin\'es par la m\'ethode la plus avantageuse. 
Le poids $P$ de cette erreur est donc \'egal \`a
$\frac{ s S p_5^{(i)2} }{ 2 S \varepsilon'^{(i)2} } $.

\setcounter{equation}{0}
\renewcommand\theequation{\Alph{equation}}

Il s'agit maintenant de d\'eterminer $S p_5 ^{(i)2}$. Pour cela, on multipliera respectivement chacune des \'equations
de condition repr\'esent\'ees par l'\'equation~(\ref{eq:orig:laplace1}),
d'abord par le coefficient du premier \'el\'ement, et l'on prendra la somme de ces produits; ensuite par le coefficient du second \'el\'ement, 
et l'on prendra la somme de ces produits, et ainsi du reste. On aura, en observant que par les conditions de la m\'ethode la plus avantageuse
$ S p^{(i)} \varepsilon^{(i)} = 0$, $S q^{(i)} \varepsilon^{(i)} = 0$, $\ldots $, les six \'equations suivantes :

\begin{equation}\label{eq:orig:laplaceA}
\left\{
\begin{array}{ccccccccccccccccccc}
\overline { p \alpha }
& = & p^{(2)} & z
& + & \overline { p q } & z'
& + & \overline { p r } & z''
& + & \overline { p t } & z'''
& + & \overline { p \gamma } & z^{iv}
& + & \overline { p \lambda } & z^{v},\\
\overline { q \alpha }
& = & \overline{ p q } & z
& + & q^{(2)} & z'
& + & \overline{ q r } & z''
& + & \overline{ q t } & z'''
& + & \overline{ q \gamma } & z^{iv}
& + & \overline{ q \lambda } & z^{v},\\
\overline { r \alpha }
& = & \overline{ r p } & z
& + & \overline{ r q } & z'
& + & r^{(2)} & z''
& + & \overline{ r t } & z'''
& + & \overline{ r \gamma } & z^{iv}
& + & \overline{ r \lambda } & z^{v},\\
\overline { t \alpha }
& = & \overline{ t p } & z
& + & \overline{ t q } & z'
& + & \overline{ t r } & z''
& + & t^{(2)} & z'''
& + & \overline{ t \gamma } & z^{iv}
& + & \overline{ t \lambda } & z^{v},\\
\overline { \gamma \alpha }
& = & \overline{ \gamma p } & z
& + & \overline{ \gamma q } & z'
& + & \overline{ \gamma r } & z''
& + & \overline{ \gamma t } & z'''
& + & \gamma^{(2)} & z^{iv}
& + & \overline{ \gamma \lambda } & z^{v},\\
\overline { \lambda \alpha }
& = & \overline{ \lambda p } & z
& + & \overline{ \lambda q } & z'
& + & \overline{ \lambda r } & z''
& + & \overline{ \lambda t } & z'''
& + & \overline{ \lambda \gamma } & z^{iv}
& + & \lambda^{(2)} & z^{v},
\end{array}
\right.
\end{equation}
o\`u l'on doit observer que nous supposons
\begin{displaymath}
p^{(2)}  = S p^{(i)2}, \quad \overline { p q }  = S p^{(i)}q^{(i)}, \quad q^{(2)}  = S q^{(i)2}, \quad \overline{qr} = S q^{(i)}r^{(i)}, \quad \ldots
\end{displaymath}
%
%
\newpage
\setcounter{equation}{2}
\renewcommand\theequation{\arabic{equation}}
\noindent
If we replace this value in Equation~(\ref{eq:orig:laplace2}) and if we perform
\begin{eqnarray}
 \nonumber  t_2^{(i)} & = & t_1^{(i)} - \gamma_1^{(i)} \frac { S \gamma_1^{(i)} t_1^{(i)} }{ S \gamma_1^{(i)2} },\\
 \nonumber  r_2^{(i)} & = & r_1^{(i)} - \gamma_1^{(i)} \frac { S \gamma_1^{(i)} r_1^{(i)} }{ S \gamma_1^{(i)2} },\\
 \nonumber  q_2^{(i)} & = & q_1^{(i)} - \gamma_1^{(i)} \frac { S \gamma_1^{(i)} q_1^{(i)} }{ S \gamma_1^{(i)2} },\\
 \nonumber  p_2^{(i)} & = & p_1^{(i)} - \gamma_1^{(i)} \frac { S \gamma_1^{(i)} p_1^{(i)} }{ S \gamma_1^{(i)2} },\\
 \nonumber  \alpha_2^{(i)} & = & \alpha_1^{(i)} - \gamma_1^{(i)} \frac { S \gamma_1^{(i)} \alpha_1^{(i)} }{ S \gamma_1^{(i)2} },
\end{eqnarray}
we have
\begin{equation}\label{eq:trad:laplace3}
\varepsilon^{(i)}
= p_2^{(i)} z
+ q_2^{(i)} z'
+ r_2^{(i)} z''
+ t_2^{(i)} z'''
- \alpha_2^{(i)}.
\end{equation}
Continuing in a similar manner, we end up with an equation of the form
\begin{equation}\label{eq:trad:laplace4}
\varepsilon^{(i)} = p_5^{(i)} z - \alpha_5^{(i)}.
\end{equation}
From n$^{\textmd{o}}$20 of \textit{Livre} II, we know that, if the value of $z$
is determined by this equation and if $u$ is the error of the value, the
probability of this error will be
\begin{displaymath}
\sqrt{ \frac{ s S p_5^{(i)2} }{ 2 S \varepsilon'^{(i)2} \pi } } 
e^{- \frac{ s S p_5^{(i)2} }{ 2 S \varepsilon'^{(i)2} } u^2 }, 
\end{displaymath}
where $ S \varepsilon'^{(i)2} $ is the sum of the squares of the residuals of the equations of condition, 
after we replaced the elements determined by \textit{la m\'ethode la plus avantageuse}. \textit{Le poids} $P$
of this error is then equal to 
$\frac{ s S p_5^{(i)2} }{ 2 S \varepsilon'^{(i)2} } $.

\setcounter{equation}{0}
\renewcommand\theequation{\Alph{equation}}

Our next task is to determine $S p_5 ^{(i)2}$. For this, we
multiply each of these equations represented by
Equation~(\ref{eq:orig:laplace1}), first by the coefficient of the first
element, and we take the sum of these products; then by the coefficient of
the second element, and we take the sum of these products, and so on for
the remaining. We have, by observing that the conditions of
\textit{la m\'ethode la plus avantageuse} $ S p^{(i)}
\varepsilon^{(i)} = 0$, $S q^{(i)} \varepsilon^{(i)} = 0$, $\ldots $, the six
following equations: 

\begin{equation}\label{eq:trad:laplaceA}
\left\{
\begin{array}{ccccccccccccccccccc}
\overline { p \alpha }
& = & p^{(2)} & z
& + & \overline { p q } & z'
& + & \overline { p r } & z''
& + & \overline { p t } & z'''
& + & \overline { p \gamma } & z^{iv}
& + & \overline { p \lambda } & z^{v},\\
\overline { q \alpha }
& = & \overline{ p q } & z
& + & q^{(2)} & z'
& + & \overline{ q r } & z''
& + & \overline{ q t } & z'''
& + & \overline{ q \gamma } & z^{iv}
& + & \overline{ q \lambda } & z^{v},\\
\overline { r \alpha }
& = & \overline{ r p } & z
& + & \overline{ r q } & z'
& + & r^{(2)} & z''
& + & \overline{ r t } & z'''
& + & \overline{ r \gamma } & z^{iv}
& + & \overline{ r \lambda } & z^{v},\\
\overline { t \alpha }
& = & \overline{ t p } & z
& + & \overline{ t q } & z'
& + & \overline{ t r } & z''
& + & t^{(2)} & z'''
& + & \overline{ t \gamma } & z^{iv}
& + & \overline{ t \lambda } & z^{v},\\
\overline { \gamma \alpha }
& = & \overline{ \gamma p } & z
& + & \overline{ \gamma q } & z'
& + & \overline{ \gamma r } & z''
& + & \overline{ \gamma t } & z'''
& + & \gamma^{(2)} & z^{iv}
& + & \overline{ \gamma \lambda } & z^{v},\\
\overline { \lambda \alpha }
& = & \overline{ \lambda p } & z
& + & \overline{ \lambda q } & z'
& + & \overline{ \lambda r } & z''
& + & \overline{ \lambda t } & z'''
& + & \overline{ \lambda \gamma } & z^{iv}
& + & \lambda^{(2)} & z^{v},
\end{array}
\right.
\end{equation}
where we have defined
\begin{displaymath}
p^{(2)}  = S p^{(i)2}, \quad \overline { p q }  = S p^{(i)}q^{(i)}, \quad q^{(2)}  = S q^{(i)2}, \quad \overline{qr} = S q^{(i)}r^{(i)}, \quad \ldots
\end{displaymath}

%
%
\newpage
\setcounter{equation}{1}
\renewcommand\theequation{\Alph{equation}}
\indent
Si l'on multiplie  pareillement les \'equations de condition repr\'esent\'ees par l'\'equation~(\ref{eq:orig:laplace2})
respectivement par les coefficients de $z$ et que l'on ajoute ces produits, ensuite par les coefficients de $z'$ en ajoutant
encore ces produits, et ainsi de suite, on aura le syst\`eme suivant d'\'equations, en observant que
$S p_1^{(i)}\varepsilon^{(i)} = 0$, $S q_1^{(i)}\varepsilon^{(i)} = 0$, $\ldots$,  par les conditions de la m\'ethode la plus avantageuse,
\begin{equation}\label{eq:orig:laplaceB}
\left\{
\begin{array}{cccccccccccccccc}
\overline { p_1 \alpha_1 }
& = & p_1^{(2)} & z
& + & \overline { p_1 q_1 } & z'
& + & \overline { p_1 r_1 } & z''
& + & \overline { p_1 t_1 } & z'''
& + & \overline { p_1 \gamma_1 } & z^{iv},\\
\overline { q_1 \alpha_1 }
& = & \overline{ p_1 q_1 } & z
& + & q_1^{(2)} & z'
& + & \overline{ q_1 r_1 } & z''
& + & \overline{ q_1 t_1 } & z'''
& + & \overline{ q_1 \gamma_1 } & z^{iv},\\
\overline { r_1 \alpha_1 }
& = & \overline{ p_1 r_1 } & z
& + & \overline{ q_1 r_1 } & z'
& + & r_1^{(2)} & z''
& + & \overline{ r_1 t_1 } & z'''
& + & \overline{ r_1 \gamma_1 } & z^{iv},\\
\overline { t_1 \alpha_1 }
& = & \overline{ p_1 t_1 } & z
& + & \overline{ q_1 t_1 } & z'
& + & \overline{ r_1 t_1 } & z''
& + & t_1^{(2)} & z'''
& + & \overline{ t_1 \gamma_1 } & z^{iv},\\
\overline { \gamma_1 \alpha_1 }
& = & \overline{ p_1 \gamma_1 } & z
& + & \overline{ q_1 \gamma_1 } & z'
& + & \overline{ r_1 \gamma_1 } & z''
& + & \overline{ t_1 \gamma_1 } & z'''
& + & \gamma_1^{(2)} & z^{iv},
\end{array}
\right.
\end{equation}
o\`u l'on doit observer que
\begin{displaymath}
 \overline { p_1 q_1 }  = S p_1^{(i)}q_1^{(i)}, \quad p_1^{(2)}  = S p_1^{(i)2}, \quad \ldots
\end{displaymath}
En substituant, au lieu de $p_1^{(i)}$, $q_1^{(i)}$, $\ldots$, leurs valeurs pr\'ec\'edentes, on a
\begin{displaymath}
\overline { p_1 q_1 }  = S p^{(i)}q^{(i)} - \frac{S\lambda^{(i)} p^{(i)} S \lambda^{(i)} q^{(i)} }{S \lambda^{(i)2} }
\end{displaymath}
ou
\begin{displaymath}
\overline { p_1 q_1 }  = \overline{ p q } - \frac{\overline{ \lambda p }\textmd{ }\overline{ \lambda q }}{\lambda^{(2)}};
\end{displaymath}
on a pareillement
\begin{displaymath}
\begin{array}{c}
\begin{array}{ccccc}
p_1^{(2)} & = & p^{(2)} & - & \frac { \overline { \lambda p }^2 }{ \lambda^{(2)} }, \\
q_1^{(2)} & = & q^{(2)} & - & \frac { \overline { \lambda q }^2 }{ \lambda^{(2)} }, \\
\overline{ p_1 r_1 } & = & \overline { p r } & - & \frac { \overline { \lambda p }\textmd{ }\overline{\lambda r} }{ \lambda^{(2)} }, \\
\end{array}\\
......................................,\\
\begin{array}{ccccc}
\overline{ p_1 \alpha_1 } & = & \overline { p \alpha } & - & \frac { \overline { \lambda p }\textmd{ }\overline{\lambda \alpha} }{ \lambda^{(2)} }, 
\end{array}\\
......................................
\end{array}
\end{displaymath}
Ainsi les coefficients du syst\`eme des \'equations~(\ref{eq:orig:laplaceB}) se d\'eduisent 
facilement des coefficients du syst\`eme des \'equations~(\ref{eq:orig:laplaceA}).

Les \'equations de condition repr\'esent\'ees par l'\'equation~(\ref{eq:orig:laplace3}) donneront semblablement
le syst\`eme suivant d'\'equations
\begin{equation}\label{eq:orig:laplaceC}
\left\{
\begin{array}{ccccccccccccc}
\overline { p_2 \alpha_2 }
& = & p_2^{(2)} & z
& + & \overline { p_2 q_2 } & z'
& + & \overline { p_2 r_2 } & z''
& + & \overline { p_2 t_2 } & z''',\\
\overline { q_2 \alpha_2 }
& = & \overline{ p_2 q_2 } & z
& + & q_2^{(2)} & z'
& + & \overline{ q_2 r_2 } & z''
& + & \overline{ q_2 t_2 } & z''',\\
\overline { r_2 \alpha_2 }
& = & \overline{ p_2 r_2 } & z
& + & \overline{ q_2 r_2 } & z'
& + & r_2^{(2)} & z''
& + & \overline{ r_2 t_2 } & z''',\\
\overline { t_2 \alpha_2 }
& = & \overline{ p_2 t_2 } & z
& + & \overline{ q_2 t_2 } & z'
& + & \overline{ r_2 t_2 } & z''
& + & t_2^{(2)} & z''',\\
\end{array}
\right.
\end{equation}
et l'on a 
\begin{displaymath}
\begin{array}{c}
\begin{array}{ccccc}
p_2^{(2)} & = & p_1^{(2)} & - & \frac { \overline { \gamma_1 p_1 }^2 }{ \gamma_1^{(2)} }, \\
\overline{ p_2 q_2 } & = & \overline { p_1 q_1 } & - & \frac { \overline { \gamma_1 p_1 }\textmd{ }\overline{q_1 \gamma_1} }{ \gamma_1^{(2)} }, \\
\end{array}\\
......................................,\\
\begin{array}{ccccc}
\overline{ p_2 \alpha_2 } & = & \overline { p_1 \alpha_1 } & - & \frac { \overline { \gamma_1 p_1 }\textmd{ }\overline{\gamma_1 \alpha_1} }{ \gamma_1^{(2)} }, \\
\end{array}\\
......................................,
\end{array}
\end{displaymath}

%
%
\newpage
\setcounter{equation}{1}
\renewcommand\theequation{\Alph{equation}}
\indent
If we multiply similarly the equations represented by Equation~(\ref{eq:trad:laplace2}) respectively by the coefficients of $z$
and we add these products, then by the coefficient of $z'$ adding again these products, and so on, we have the following 
system of equations, by noting that 
$S p_1^{(i)}\varepsilon^{(i)} = 0$, $S q_1^{(i)}\varepsilon^{(i)} = 0$, $\ldots$, from the conditions of \textit{la m\'ethode la plus avantageuse},
\begin{equation}\label{eq:trad:laplaceB}
\left\{
\begin{array}{cccccccccccccccc}
\overline { p_1 \alpha_1 }
& = & p_1^{(2)} & z
& + & \overline { p_1 q_1 } & z'
& + & \overline { p_1 r_1 } & z''
& + & \overline { p_1 t_1 } & z'''
& + & \overline { p_1 \gamma_1 } & z^{iv},\\
\overline { q_1 \alpha_1 }
& = & \overline{ p_1 q_1 } & z
& + & q_1^{(2)} & z'
& + & \overline{ q_1 r_1 } & z''
& + & \overline{ q_1 t_1 } & z'''
& + & \overline{ q_1 \gamma_1 } & z^{iv},\\
\overline { r_1 \alpha_1 }
& = & \overline{ p_1 r_1 } & z
& + & \overline{ q_1 r_1 } & z'
& + & r_1^{(2)} & z''
& + & \overline{ r_1 t_1 } & z'''
& + & \overline{ r_1 \gamma_1 } & z^{iv},\\
\overline { t_1 \alpha_1 }
& = & \overline{ p_1 t_1 } & z
& + & \overline{ q_1 t_1 } & z'
& + & \overline{ r_1 t_1 } & z''
& + & t_1^{(2)} & z'''
& + & \overline{ t_1 \gamma_1 } & z^{iv},\\
\overline { \gamma_1 \alpha_1 }
& = & \overline{ p_1 \gamma_1 } & z
& + & \overline{ q_1 \gamma_1 } & z'
& + & \overline{ r_1 \gamma_1 } & z''
& + & \overline{ t_1 \gamma_1 } & z'''
& + & \gamma_1^{(2)} & z^{iv},
\end{array}
\right.
\end{equation}
where we have defined
\begin{displaymath}
 \overline { p_1 q_1 }  = S p_1^{(i)}q_1^{(i)}, \quad p_1^{(2)}  = S p_1^{(i)2}, \quad \ldots
\end{displaymath}
Substituting $p_1^{(i)}$, $q_1^{(i)}$, $\ldots$ with their previous values, we have
\begin{displaymath}
\overline { p_1 q_1 }  = S p^{(i)}q^{(i)} - \frac{S\lambda^{(i)} p^{(i)} S \lambda^{(i)} q^{(i)} }{S \lambda^{(i)2} }
\end{displaymath}
or 
\begin{displaymath}
\overline { p_1 q_1 }  = \overline{ p q } - \frac{\overline{ \lambda p }\textmd{ }\overline{ \lambda q }}{\lambda^{(2)}};
\end{displaymath}
we have similarly
\begin{displaymath}
\begin{array}{c}
\begin{array}{ccccc}
p_1^{(2)} & = & p^{(2)} & - & \frac { \overline { \lambda p }^2 }{ \lambda^{(2)} }, \\
q_1^{(2)} & = & q^{(2)} & - & \frac { \overline { \lambda q }^2 }{ \lambda^{(2)} }, \\
\overline{ p_1 r_1 } & = & \overline { p r } & - & \frac { \overline { \lambda p }\textmd{ }\overline{\lambda r} }{ \lambda^{(2)} }, \\
\end{array}\\
......................................,\\
\begin{array}{ccccc}
\overline{ p_1 \alpha_1 } & = & \overline { p \alpha } & - & \frac { \overline { \lambda p }\textmd{ }\overline{\lambda \alpha} }{ \lambda^{(2)} }, 
\end{array}\\
......................................
\end{array}
\end{displaymath}
Doing so, the coefficients of the system of equations~(\ref{eq:trad:laplaceB}) are easily
computable from the coefficients of the system of equations~(\ref{eq:trad:laplaceA}).

The equations represented by Equation~(\ref{eq:trad:laplace3}) similarly give the following system of equations
\begin{equation}\label{eq:trad:laplaceC}
\left\{
\begin{array}{ccccccccccccc}
\overline { p_2 \alpha_2 }
& = & p_2^{(2)} & z
& + & \overline { p_2 q_2 } & z'
& + & \overline { p_2 r_2 } & z''
& + & \overline { p_2 t_2 } & z''',\\
\overline { q_2 \alpha_2 }
& = & \overline{ p_2 q_2 } & z
& + & q_2^{(2)} & z'
& + & \overline{ q_2 r_2 } & z''
& + & \overline{ q_2 t_2 } & z''',\\
\overline { r_2 \alpha_2 }
& = & \overline{ p_2 r_2 } & z
& + & \overline{ q_2 r_2 } & z'
& + & r_2^{(2)} & z''
& + & \overline{ r_2 t_2 } & z''',\\
\overline { t_2 \alpha_2 }
& = & \overline{ p_2 t_2 } & z
& + & \overline{ q_2 t_2 } & z'
& + & \overline{ r_2 t_2 } & z''
& + & t_2^{(2)} & z''',\\
\end{array}
\right.
\end{equation}
and we have
\begin{displaymath}
\begin{array}{c}
\begin{array}{ccccc}
p_2^{(2)} & = & p_1^{(2)} & - & \frac { \overline { \gamma_1 p_1 }^2 }{ \gamma_1^{(2)} }, \\
\overline{ p_2 q_2 } & = & \overline { p_1 q_1 } & - & \frac { \overline { \gamma_1 p_1 }\textmd{ }\overline{q_1 \gamma_1} }{ \gamma_1^{(2)} }, \\
\end{array}\\
......................................,\\
\begin{array}{ccccc}
\overline{ p_2 \alpha_2 } & = & \overline { p_1 \alpha_1 } & - & \frac { \overline { \gamma_1 p_1 }\textmd{ }\overline{\gamma_1 \alpha_1} }{ \gamma_1^{(2)} }, \\
\end{array}\\
......................................,
\end{array}
\end{displaymath}

%
%
\newpage
\setcounter{equation}{3}
\renewcommand\theequation{\Alph{equation}}
\indent
On a pareillement le syst\`eme d'\'equations
\begin{equation}\label{eq:trad:laplaceD}
\left\{
\begin{array}{cccccccccc}
\overline { p_3 \alpha_3 }
& = & p_3^{(2)} & z
& + & \overline { p_3 q_3 } & z'
& + & \overline { p_3 r_3 } & z'',\\
\overline { q_3 \alpha_3 }
& = & \overline{ p_3 q_3 } & z
& + & q_3^{(2)} & z'
& + & \overline{ q_3 r_3 } & z'',\\
\overline { r_3 \alpha_3 }
& = & \overline{ p_3 r_3 } & z
& + & \overline{ q_3 r_3 } & z'
& + & r_3^{(2)} & z'',\\
\end{array}
\right.
\end{equation}
en faisant
\begin{displaymath}
\begin{array}{c}
\begin{array}{ccccc}
p_3^{(2)} & = & p_2^{(2)} & - & \frac { \overline { p_2 t_2 }^2 }{ t_2^{(2)} }, \\
\overline{ p_3 q_3 } & = & \overline { p_2 q_2 } & - & \frac { \overline { p_2 t_2 }\textmd{ }\overline{q_2 t_2} }{ t_2^{(2)} }, \\
\overline{ p_3 \alpha_3 } & = & \overline { p_2 \alpha_2 } & - & \frac { \overline { t_2 p_2 }\textmd{ }\overline{t_2 \alpha_2} }{ t_2^{(2)} }, \\
\end{array}\\
......................................;
\end{array}
\end{displaymath}
on aura encore
\begin{equation}\label{eq:trad:laplaceE}
\left\{
\begin{array}{ccccccc}
\overline { p_4 \alpha_4 }
& = & p_4^{(2)} & z
& + & \overline { p_4 q_4 } & z', \\
\overline { q_4 \alpha_4 }
& = & \overline{ p_4 q_4 } & z
& + & q_4^{(2)} & z', \\
\end{array}
\right.
\end{equation}
en faisant
\begin{displaymath}
\begin{array}{c}
\begin{array}{ccccc}
p_4^{(2)} & = & p_3^{(2)} & - & \frac { \overline { p_3 r_3 }^2 }{ r_3^{(2)} }, \\
\overline{ p_4 q_4 } & = & \overline { p_3 q_3 } & - & \frac { \overline { p_3 r_3 }\textmd{ }\overline{q_3 r_3} }{ r_3^{(2)} }, \\
\overline{ p_4 \alpha_4 } & = & \overline { p_3 \alpha_3 } & - & \frac { \overline { p_3 r_3 }\textmd{ }\overline{ \alpha_3 r_3 } }{ r_3^{(2)} }, \\
\end{array}\\
......................................
\end{array}
\end{displaymath}
Enfin on aura
\begin{equation}\label{eq:trad:laplaceF}
\overline { p_5 \alpha_5 } = p_5^{(2)} z, 
\end{equation}
en faisant
\begin{displaymath}
p_5^{(2)} = p_4^{(2)} - \frac { \overline { p_4 q_4 }^2 }{ q_4^{(2)} }, \quad
\overline{ p_5 \alpha_5 } = \overline { p_4 \alpha_4 } - \frac { \overline { p_4 q_4 }\textmd{ }\overline{ q_4 \alpha_4 } }{ q_4^{(2)} };
\end{displaymath}
$p_5^{(2)}$ est la valeur $S p_5^{(i)2}$, et le poids $P$ sera
\begin{displaymath}
\frac{sp_5^{(2)}}{2S\varepsilon'^{(i)2}}.
\end{displaymath}
On voit par la suite des valeurs $p^{(2)}$, $p_1^{(2)}$, $p_2^{(2)}$, $\ldots$
qu'elles vont en diminuant sans cesse, et qu'ainsi, pour le m\^eme nombre
d'observations, le poids $P$ diminue quand le nombre des \'el\'ements augmente. 

Si l'on consid\`ere la suite des \'equations qui d\'eterminent 
$\overline{p_5 \alpha_5}$, on voit que cette fonction, d\'evelopp\'ee suivant les coefficients du syst\`eme des \'equations~(\ref{eq:trad:laplaceA}),
est de la forme
\begin{displaymath}
\overline{p \alpha} + M \overline{q \alpha} + N \overline{r \alpha} + \ldots,
\end{displaymath}
le coefficient de $\overline{p \alpha}$ \'etant l'unit\'e. Il suit de l\`a que si l'on r\'esout les \'equations~(\ref{eq:trad:laplaceA}),
en y laissant  $\overline{p \alpha}$, $\overline{q \alpha}$, $ \overline{r \alpha}$, $\ldots$
comme ind\'etermin\'ees, $\frac{1}{p_5^{(2)}}$ sera, en vertu de l'\'equation~(\ref{eq:trad:laplaceF}),
le coefficient de $\overline{p \alpha}$ dans l'expression de $z$. 
Pareillement, $\frac{1}{q_5^{(2)}}$ sera le coefficient de $\overline{q \alpha}$ dans l'expression de $z'$; 
$\frac{1}{r_5^{(2)}}$ sera le coefficient de $\overline{r \alpha}$ dans l'expression de $z''$; et ainsi de suite du reste;
ce qui donne un moyen de simple d'obtenir $p_5^{(2)}$, $q_5^{(2)}$, $\ldots$; mais il est plus simple encore de les determiner ainsi.

%
%
\newpage
\setcounter{equation}{3}
\renewcommand\theequation{\Alph{equation}}
\indent
We similarly have the system of equations
\begin{equation}\label{eq:orig:laplaceD}
\left\{
\begin{array}{cccccccccc}
\overline { p_3 \alpha_3 }
& = & p_3^{(2)} & z
& + & \overline { p_3 q_3 } & z'
& + & \overline { p_3 r_3 } & z'',\\
\overline { q_3 \alpha_3 }
& = & \overline{ p_3 q_3 } & z
& + & q_3^{(2)} & z'
& + & \overline{ q_3 r_3 } & z'',\\
\overline { r_3 \alpha_3 }
& = & \overline{ p_3 r_3 } & z
& + & \overline{ q_3 r_3 } & z'
& + & r_3^{(2)} & z'',\\
\end{array}
\right.
\end{equation}
by doing
\begin{displaymath}
\begin{array}{c}
\begin{array}{ccccc}
p_3^{(2)} & = & p_2^{(2)} & - & \frac { \overline { p_2 t_2 }^2 }{ t_2^{(2)} }, \\
\overline{ p_3 q_3 } & = & \overline { p_2 q_2 } & - & \frac { \overline { p_2 t_2 }\textmd{ }\overline{q_2 t_2} }{ t_2^{(2)} }, \\
\overline{ p_3 \alpha_3 } & = & \overline { p_2 \alpha_2 } & - & \frac { \overline { t_2 p_2 }\textmd{ }\overline{t_2 \alpha_2} }{ t_2^{(2)} }, \\
\end{array}\\
......................................;
\end{array}
\end{displaymath}
we also have
\begin{equation}\label{eq:orig:laplaceE}
\left\{
\begin{array}{ccccccc}
\overline { p_4 \alpha_4 }
& = & p_4^{(2)} & z
& + & \overline { p_4 q_4 } & z', \\
\overline { q_4 \alpha_4 }
& = & \overline{ p_4 q_4 } & z
& + & q_4^{(2)} & z', \\
\end{array}
\right.
\end{equation}
by doing
\begin{displaymath}
\begin{array}{c}
\begin{array}{ccccc}
p_4^{(2)} & = & p_3^{(2)} & - & \frac { \overline { p_3 r_3 }^2 }{ r_3^{(2)} }, \\
\overline{ p_4 q_4 } & = & \overline { p_3 q_3 } & - & \frac { \overline { p_3 r_3 }\textmd{ }\overline{q_3 r_3} }{ r_3^{(2)} }, \\
\overline{ p_4 \alpha_4 } & = & \overline { p_3 \alpha_3 } & - & \frac { \overline { p_3 r_3 }\textmd{ }\overline{ \alpha_3 r_3 } }{ r_3^{(2)} }, \\
\end{array}\\
......................................
\end{array}
\end{displaymath}
Finally we have
\begin{equation}\label{eq:orig:laplaceF}
\overline { p_5 \alpha_5 } = p_5^{(2)} z, 
\end{equation}
by doing
\begin{displaymath}
p_5^{(2)} = p_4^{(2)} - \frac { \overline { p_4 q_4 }^2 }{ q_4^{(2)} }, \quad
\overline{ p_5 \alpha_5 } = \overline { p_4 \alpha_4 } - \frac { \overline { p_4 q_4 }\textmd{ }\overline{ q_4 \alpha_4 } }{ q_4^{(2)} };
\end{displaymath}
$p_5^{(2)}$ is the value $S p_5^{(i)2}$, and \textit{le poids} $P$ is
\begin{displaymath}
\frac{sp_5^{(2)}}{2S\varepsilon'^{(i)2}}.
\end{displaymath}
We see from the sequence of values $p^{(2)}$, $p_1^{(2)}$, $p_2^{(2)}$, $\ldots$
that they always go diminishing, and so, for the same number of observations, 
\textit{le poids} $P$ decreases when the number of elements increases.

If we consider the sequence of equations which determine
$\overline{p_5 \alpha_5}$, we see that this function, developed according to the coefficients of the system of equations~(\ref{eq:orig:laplaceA}),
is of the form
\begin{displaymath}
\overline{p \alpha} + M \overline{q \alpha} + N \overline{r \alpha} + \ldots,
\end{displaymath}
the coefficient of $\overline{p \alpha}$ being the unity. It follows from there that if we solve
the equations~(\ref{eq:orig:laplaceA}), by leaving
$\overline{p \alpha}$, $\overline{q \alpha}$, $ \overline{r \alpha}$, $\ldots$ as unknowns, 
$\frac{1}{p_5^{(2)}}$ is, due to Equation~(\ref{eq:orig:laplaceF}),
the coefficient of $\overline{p \alpha}$ in the expression of $z$. 
Similarly, $\frac{1}{q_5^{(2)}}$ is the coefficient of $\overline{q \alpha}$ in the expression of $z'$; 
$\frac{1}{r_5^{(2)}}$ is the coefficient of $\overline{r \alpha}$ in the expression of $z''$; and so on for the others;
this gives a simple mean to obtain $p_5^{(2)}$, $q_5^{(2)}$, $\ldots$; but it is even simpler to compute them as follows.

%
%
\newpage
\setcounter{equation}{6}
\renewcommand\theequation{\Alph{equation}}
\indent

D'abord l'\'equation~(\ref{eq:orig:laplaceF}) donne la valeur de $p_5^{(2)}$ et
de $z$. Si dans le syst\`eme des \'equations~(\ref{eq:orig:laplaceE}) on \'elimine
$z$ au lieu de $z'$, on aura une seule \'equation en $z'$, de la forme
\begin{displaymath}
\overline { q_5 \alpha_5 }  =  q_5^{(2)}  z';
\end{displaymath}
en faisant
\begin{displaymath}
q_5^{(2)} = q_4^{(2)} - \frac { \overline { p_4 q_4 }^2 }{ p_4^{(2)} }, \quad
\overline{ q_5 \alpha_5 } = \overline { q_4 \alpha_4 } - \frac { \overline { p_4 q_4 }\textmd{ }\overline{ p_4 \alpha_4 } }{ p_4^{(2)} }.
\end{displaymath}

Si dans le syst\`eme des \'equations~(\ref{eq:orig:laplaceD}) on \'elimine $z$
au lieu de $z''$, pour ne conserver \`a la fin du calcul que $z''$, on aura
$r_5^{(2)}$ en changeant dans la suite des \'equations qui, \`a partir de ce
syst\`eme, d\'eterminent $p_5^{(2)}$, la lettre $p$ dans la lettre $r$, et
r\'eciproquement. On aura ainsi
\begin{displaymath}
\begin{array}{c}
\begin{array}{ccccc}
r_4^{(2)} & = & r_3^{(2)} & - & \frac { \overline { p_3 r_3 }^2 }{ p_3^{(2)} }, \\
\overline{ r_4 q_4 } & = & \overline { r_3 q_3 } & - & \frac { \overline { p_3 q_3 }\textmd{ }\overline{p_3 r_3} }{ p_3^{(2)} }, \\
q_4^{(2)} & = & q_3^{(2)} & - & \frac { \overline { p_3 q_3 }^2 }{ p_3^{(2)} }, \\
r_5^{(2)} & = & r_4^{(2)} & - & \frac { \overline { p_4 q_4 }^2 }{ q_4^{(2)} }, \\
\end{array}\\
......................................
\end{array}
\end{displaymath}
Pour avoir $t_5^{(2)}$, on partira du syst\`eme~(\ref{eq:orig:laplaceC}), en
changeant, dans la suite des valeurs de $p_3^{(2)}$, $\overline{p_3q_3}$,
$\ldots$, $r_3^{(2)}$, $\overline{q_3r_3}$, $\ldots$, la lettre $p$ dans la
lettre $t$, et r\'eciproquement.

On aura pareillement la valeur de $\gamma_5^{(2)}$, en partant du syst\`eme des
\'equations~(\ref{eq:orig:laplaceB}) et changeant dans la suite des valeurs de
$p_2^{(2)}$, $p_3^{(2)}$, $\ldots$, la lettre $p$ dans la lettre $\gamma$, et
r\'eciproquement.

Enfin, on aura la valeur de $\lambda_5^{(2)}$ en changeant, dans la suite des
valeurs de $p_1^{(2)}$, $p_2^{(2)}$, $\ldots$, la lettre $p$ dans la lettre
$\lambda$, et r\'eciproquement.

%
%
\newpage
\setcounter{equation}{6}
\renewcommand\theequation{\Alph{equation}}
\indent

Firstly Equation~(\ref{eq:trad:laplaceF}) gives the value of $p_5^{(2)}$ and of $z$.
If in the system of equations~(\ref{eq:trad:laplaceE}) we eliminate 
$z$ instead of $z'$, we have a single equation in $z'$, of the form
\begin{displaymath}
\overline { q_5 \alpha_5 }  =  q_5^{(2)}  z';
\end{displaymath}
by doing
\begin{displaymath}
q_5^{(2)} = q_4^{(2)} - \frac { \overline { p_4 q_4 }^2 }{ p_4^{(2)} }, \quad
\overline{ q_5 \alpha_5 } = \overline { q_4 \alpha_4 } - \frac { \overline { p_4 q_4 }\textmd{ }\overline{ p_4 \alpha_4 } }{ p_4^{(2)} }.
\end{displaymath}

If in the system of equations~(\ref{eq:trad:laplaceD}) we eliminate 
$z$ instead of $z''$, in order to only keep at the end of the computation $z''$, we have
$r_5^{(2)}$ by changing in the sequence of equations which, starting from this system, determine
 $p_5^{(2)}$, the letter $p$ by the letter $r$, and reciprocally. We then have
\begin{displaymath}
\begin{array}{c}
\begin{array}{ccccc}
r_4^{(2)} & = & r_3^{(2)} & - & \frac { \overline { p_3 r_3 }^2 }{ p_3^{(2)} }, \\
\overline{ r_4 q_4 } & = & \overline { r_3 q_3 } & - & \frac { \overline { p_3 q_3 }\textmd{ }\overline{p_3 r_3} }{ p_3^{(2)} }, \\
q_4^{(2)} & = & q_3^{(2)} & - & \frac { \overline { p_3 q_3 }^2 }{ p_3^{(2)} }, \\
r_5^{(2)} & = & r_4^{(2)} & - & \frac { \overline { p_4 q_4 }^2 }{ q_4^{(2)} }, \\
\end{array}\\
......................................
\end{array}
\end{displaymath}
In order to have $t_5^{(2)}$, we start from the
system~(\ref{eq:trad:laplaceC}), by changing, in the sequence of values of
$p_3^{(2)}$, $\overline{p_3q_3}$, $\ldots$, $r_3^{(2)}$, $\overline{q_3r_3}$,
$\ldots$, the letter $p$ by the letter $t$, and reciprocally.

We similarly have the value of $\gamma_5^{(2)}$, starting from the system
of equations~(\ref{eq:trad:laplaceB}) and changing in the sequence of values of
$p_2^{(2)}$, $p_3^{(2)}$, $\ldots$, the letter $p$ by the letter $\gamma$, and
reciprocally.

Finally, we have the value of $\lambda_5^{(2)}$ by changing, in the sequence
of values of $p_1^{(2)}$, $p_2^{(2)}$, $\ldots$, the letter $p$ by the letter
$\lambda$, and reciprocally.

\newpage

%
%
\newpage
\indent

5. Appliquons maintenant cette m\'ethode \`a un exemple. Pour cela, j'ai
profit\'e de l'immense travail que Bouvard vient de terminer sur les mouvements
de Jupiter et de Saturne, dont il a construit des Tables tr\`es pr\'ecises. Il a
fait usage de toutes les oppositions observ\'ees par Bradley et par les
astronomes qui l'ont suivi~: il les a discut\'ees de nouveau et avec le plus
gran soin, ce qui lui a donn\'e $126$ \'equations de condition pour le
mouvement de Jupiter en longitude et $129$ \'equations pour le mouvement de
Saturne. Dans ces derni\`eres \'equations, Bouvard a fait entrer la masse
d'Uranus comme ind\'etermin\'ee. Voici les \'equations finales qu'il a conclues
par la m\'ethode la plus avantageuse :
\begin{displaymath}
\begin{array}{rcl}
    7212",600 & = &   795938 z -  12729398 z' +   6788,2 z'' -  1959,0 z''' + 696,13 z^{iv} +   2602 z^{v}, \\
 -738297",800 & = &-12729398 z + 424865729 z' - 153106,5 z'' - 39749,1 z''' -   5459 z^{iv} +   5722 z^{v}, \\
     237",782 & = &   6788,2 z -  153106,5 z' +  71,8720 z'' -  3,2252 z''' + 1,2484 z^{iv} + 1,3371 z^{v}, \\
     -40",335 & = &  -1959,0 z -   39749,1 z' -   3,2252 z'' + 57,1911 z''' + 3,6213 z^{iv} + 1,1128 z^{v}, \\
    -343",455 & = &   696,13 z -      5459 z' +   1,2484 z'' +  3,6213 z''' + 21,543 z^{iv} + 46,310 z^{v}, \\
   -1002",900 & = &     2602 z +      5722 z' +   1,3371 z'' +  1,1128 z''' + 46,310 z^{iv} +    129 z^{v}. \\
\end{array}	  
\end{displaymath}

Dans ces equations, la masse d'Uranus est suppos\'ee $\frac{1+z}{19504}$; la
masse de Jupiter est suppos\'ee $\frac{1+z'}{1067,09}$; $z''$ est le produit de
l'\'equation du centre par la correction du p\'erih\'elie employ\'e d'abord par
Bouvard; $z'''$ est la correction de l'\'equation du centre; $z^{iv}$ est la
correction s\'eculaire du moyen mouvement; $z^v$ est la correction de
l'\'epoque de la longitude au commencement de $1750$. La seconde du degr\'e
d\'ecimal est prise pour unit\'e.

Au moyen des \'equations pr\'ec\'edentes renferm\'ees dans le syst\`eme~(\ref{eq:orig:laplaceA}), 
j'ai conclu les suivantes, renferm\'ees dans le syst\`eme~(\ref{eq:orig:laplaceB}) :
\begin{displaymath}
\begin{array}{rcl}
    27441",68 & = &     743454 z -    12844814 z' +     6761,23 z'' -    1981,45 z''' -    237,97 z^{iv},  \\
  -693812",58 & = &  -12844814 z +   424611920 z' -   153165,81 z'' -   39798,46 z''' -   7513,15 z^{iv},  \\
    248",1772 & = &    6761,23 z -   153165,81 z' +     71,8581 z'' -     3,2367 z''' +    0,7684 z^{iv},  \\
    -31",6836 & = &   -1981,45 z -    39798,46 z' -      3,2367 z'' +    57,1815 z''' +    3,2218 z^{iv},  \\
     16",5783 & = &    -237,97 z -     7513,15 z' +      0,7684 z'' +     3,2218 z''' +    4,9181 z^{iv} . \\
\end{array}	  
\end{displaymath}

De ces \'equations, j'ai tir\'e les quatre suivantes, renferm\'ees dans le syst\`eme~(\ref{eq:orig:laplaceC}),
\begin{displaymath}
\begin{array}{rcl}
    28243",85  & = &   731939,5 z -   13208350 z' +    6798,41 z'' -   1825,56 z''' , \\
  -668486",70  & = &  -13208350 z +  413134432 z' -   151992,0 z'' -   34876,7 z''' , \\
    245",5870  & = &    6798,41 z -   151992,0 z' +    71,7381 z'' -    3,7401 z''' , \\
     42",5434  & = &   -1825,56 z -    34876,7 z' -     3,7401 z'' +   55,0710 z''' ; \\
\end{array}	  
\end{displaymath}
ces derni\`eres \'equations donnent les suivantes, renferm\'ees dans le syst\`eme~(\ref{eq:orig:laplaceD}),
\begin{displaymath}
\begin{array}{rcl}
   26833",55 & = &  671414,7 z -  14364541 z' +   6674,43 z'' , \\
  -695430",0 & = & -14364541 z + 391046861 z' -  154360,6 z'' , \\
   242",6977 & = &   6674,43 z -  154360,6 z' +   71,4841 z'' . \\
\end{array}	  
\end{displaymath}
Enfin j'ai conclu de l\`a les deux \'equations,  renferm\'ees dans le syst\`eme~(\ref{eq:orig:laplaceE}) :
\begin{displaymath}
	  4172",95  =  48442 z +    48020 z', \quad -171455",2  =  48020 z + 57725227 z' .
\end{displaymath}

%
%
\newpage
\indent

5. We now apply this method to an example. For this, I have benefited from
the immense work that Bouvard has just finished on the movements of Jupiter and
Saturn, from which he has constructed extremely accurate Tables.  He has
used all the observations from Bradley and from the astronomers that have
followed him: he has discussed them again and with the greatest care, which has
given him $126$ equations for the movement of Jupiter in longitude and $129$
equations for the movement of Saturn. In these equations, Bouvard has
introduced the mass of Uranus as unknown.  Here are the final equations
that he has obtained by \textit{la m\'ethode la plus avantageuse}:

\begin{displaymath}
\begin{array}{rcl}
    7212".600 & = &   795938 z -  12729398 z' +   6788.2 z'' -  1959.0 z''' + 696.13 z^{iv} +   2602 z^{v}, \\
 -738297".800 & = &-12729398 z + 424865729 z' - 153106.5 z'' - 39749.1 z''' -   5459 z^{iv} +   5722 z^{v}, \\
     237".782 & = &   6788.2 z -  153106.5 z' +  71.8720 z'' -  3.2252 z''' + 1.2484 z^{iv} + 1.3371 z^{v}, \\
     -40".335 & = &  -1959.0 z -   39749.1 z' -   3.2252 z'' + 57.1911 z''' + 3.6213 z^{iv} + 1.1128 z^{v}, \\
    -343".455 & = &   696.13 z -      5459 z' +   1.2484 z'' +  3.6213 z''' + 21.543 z^{iv} + 46.310 z^{v}, \\
   -1002".900 & = &     2602 z +      5722 z' +   1.3371 z'' +  1.1128 z''' + 46.310 z^{iv} +    129 z^{v}. \\
\end{array}	  
\end{displaymath}

In these equations, the mass of Uranus is supposed to be $\frac{1+z}{19504}$;
the mass of Jupiter is supposed to be $\frac{1+z'}{1067.09}$; $z''$ is the
product if the equation of the center by the correction of the periapsis
firstly employed by Bouvard; $z'''$ is the correction of the equation of the
center; $z^{iv}$ is the secular correction of the mean movement; $z^{v}$ is the
correction of the epoch of the longitude beginning in $1750$. The second of the
decimal degree is taken for unit.

With these former equations corresponding to the system~(\ref{eq:trad:laplaceA}), 
I have concluded the followings, corresponding to the system~(\ref{eq:trad:laplaceB}) :
\begin{displaymath}
\begin{array}{rcl}
    27441".68 & = &     743454 z -    12844814 z' +     6761.23 z'' -    1981.45 z''' -    237.97 z^{iv},  \\
  -693812".58 & = &  -12844814 z +   424611920 z' -   153165.81 z'' -   39798.46 z''' -   7513.15 z^{iv},  \\
    248".1772 & = &    6761.23 z -   153165.81 z' +     71.8581 z'' -     3.2367 z''' +    0.7684 z^{iv},  \\
    -31".6836 & = &   -1981.45 z -    39798.46 z' -      3.2367 z'' +    57.1815 z''' +    3.2218 z^{iv},  \\
     16".5783 & = &    -237.97 z -     7513.15 z' +      0.7684 z'' +     3.2218 z''' +    4.9181 z^{iv} . \\
\end{array}	  
\end{displaymath}

From these equations, I have drawn the four followings, corresponding to the system~(\ref{eq:trad:laplaceC}),
\begin{displaymath}
\begin{array}{rcl}
    28243".85  & = &   731939.5 z -   13208350 z' +    6798.41 z'' -   1825.56 z''' , \\
  -668486".70  & = &  -13208350 z +  413134432 z' -   151992.0 z'' -   34876.7 z''' , \\
    245".5870  & = &    6798.41 z -   151992.0 z' +    71.7381 z'' -    3.7401 z''' , \\
     42".5434  & = &   -1825.56 z -    34876.7 z' -     3.7401 z'' +   55.0710 z''' ; \\
\end{array}	  
\end{displaymath}
these latest equations give us the followings, corresponding to the system~(\ref{eq:trad:laplaceD}),
\begin{displaymath}
\begin{array}{rcl}
   26833".55 & = &  671414.7 z -  14364541 z' +   6674.43 z'' , \\
  -695430".0 & = & -14364541 z + 391046861 z' -  154360.6 z'' , \\
   242".6977 & = &   6674.43 z -  154360.6 z' +   71.4841 z'' . \\
\end{array}	  
\end{displaymath}
Finally I have concluded from there the two equations, corresponding to the system~(\ref{eq:trad:laplaceE}):
\begin{displaymath}
	  4172".95  =  48442 z +    48020 z', \quad -171455".2  =  48020 z + 57725227 z' .
\end{displaymath}

%
%
\newpage
\indent

Je m'arr\^ete \`a ce syst\`eme, parce qu'il est facile d'en conclude les valeurs du poids $P$ relatives aux deux \'el\'ements
$z$ et $z'$ que je d\'esirais particuli\`erement de conna\^{\i}tre. Les formules du n$^\textmd{o}$ 3 donnent, pour $z$,
\begin{displaymath}
P = \frac{s}{2S\varepsilon'^{(i)2}}\left[ 48442 - \frac{ (48020 )^2 }{57725227} \right]
\end{displaymath}
et, pour $z'$,
\begin{displaymath}
P = \frac{s}{2S\varepsilon'^{(i)2}}\left[ 57725227 - \frac{ (48020 )^2 }{ 48442 } \right].
\end{displaymath}
Le nombre $s$ des observations est ici $129$ et Bouvard a trouv\'e
\begin{displaymath}
S\varepsilon'^{(i)2} = 31096;
\end{displaymath}
on a donc, pour $z$,
\begin{displaymath}
\log  P = 2,0013595;
\end{displaymath}
et, pour $z'$,
\begin{displaymath}
\log  P = 5,0778624.
\end{displaymath}
Les \'equations pr\'ec\'edentes donnent
\begin{eqnarray}
\nonumber &z' = -0,00305,&\\
\nonumber &z = 0,08916.&
\end{eqnarray}

La masse de Jupiter est $\frac{1}{1067,09}(1+z')$. En substituant pour $z'$ sa
valeur p\'ec\'edente, cette masse devient $\frac{1}{1070,35}$. La masse du
Soleil est prise pour unit\'e. La probabilit\'e que l'erreur de $z'$ est
comprise dans les limites $\pm U$ est, par le n$^\textmd{o}$ 1,
\begin{displaymath}
\frac{\sqrt{P}}{\sqrt{\pi}} \int du\textmd{ } e^{-Pu^2},
\end{displaymath}
l'int\'egrale \'etant prise depuis $u=-U$ jusqu'a $u=U$. On trouve ainsi la
probabilit\'e que la masse de Jupiter est comprise dans les limites
\begin{displaymath}
\frac{1}{1070,35} \pm \frac{1}{100}\frac{1}{1067,09},
\end{displaymath}
\'egale \`a $\frac{1000000}{1000001}$; en sorte qu'il y a un million \`a tr\`es
peu \`a parier contre un que la valeur $\frac{1}{1070,35}$ n'est pas en erreur
d'un centi\`eme de sa valeur; ou, ce qui revient \`a fort peu pr\`es au m\^eme,
qu'apr\`es un si\`ecle de nouvelles observations, ajout\'ees aux
pr\'ec\'edentes et dicut\'ees de la m\^eme mani\`ere, le nouveau r\'esultat ne
diff\'erera pas du pr\'ec\'edent d'un centi\`eme de sa valeur.

Newton avait trouv\'e, par les observations de Pound, sur les \'elongations des
satellites de Jupiter, la masse de cette plan\`ete \' egale \`a la $1067^\textmd{e}$
partie de celle du Soleil, ce qui diff\`ere tr\`es peu du r\'esultat de Bouvard.

La masse d'Uranus est $\frac{1+z}{19504}$. En substituant pour $z$ sa
valeur p\'ec\'edente, cette masse devient $\frac{1}{17907}$. 
La probabilit\'e que cette valeur est
comprise dans les limites 
\begin{displaymath}
\frac{1}{17907} \pm \frac{1}{4}\frac{1}{19504}
\end{displaymath}
est \'egale \`a $\frac{2508}{2509}$, et la probabilit\'e que cette masse est comprise dans les limites
\begin{displaymath}
\frac{1}{17907} \pm \frac{1}{5}\frac{1}{19504}
\end{displaymath}
est \'egale \`a $\frac{215,6}{216,6}$.

%
%
\newpage
\indent
\label{hhhhere}

I stop with this system, because it is easy to conclude from it the values of the \textit{poids} $P$ corresponding
to the two elements $z$ and $z'$ which I particularly wish to know. The formula from n$^\textmd{o}$ 3 give, for $z$,
\begin{displaymath}
P = \frac{s}{2S\varepsilon'^{(i)2}}\left[ 48442 - \frac{ (48020 )^2 }{57725227} \right]
\end{displaymath}
and, for $z'$,
\begin{displaymath}
P = \frac{s}{2S\varepsilon'^{(i)2}}\left[ 57725227 - \frac{ (48020 )^2 }{ 48442 } \right].
\end{displaymath}
The number $s$ of observations is here $129$ and Bouvard has found
\begin{displaymath}
S\varepsilon'^{(i)2} = 31096;
\end{displaymath}
we then have, for $z$,
\begin{displaymath}
\log  P = 2.0013595;
\end{displaymath}
and, for $z'$,
\begin{displaymath}
\log  P = 5.0778624.
\end{displaymath}
The former equations give
\begin{eqnarray}
\nonumber &z' = -0.00305,&\\
\nonumber &z = 0.08916.&
\end{eqnarray}

The mass of Jupiter is $\frac{1}{1067.09}(1+z')$. Replacing $z'$ by its former value, 
this mass becomes $\frac{1}{1070.35}$. The mass of the Sun is taken as unity. The probability
that the error in $z'$ is between the limit $\pm U$ is, from n$^\textmd{o}$ 1,
\begin{displaymath}
\frac{\sqrt{P}}{\sqrt{\pi}} \int du\textmd{ } e^{-Pu^2},
\end{displaymath}
the integral being taken from $u=-U$ to $u=U$. We then find that the probability
for the mass of Jupiter to be between the limits
\begin{displaymath}
\frac{1}{1070.35} \pm \frac{1}{100}\frac{1}{1067.09},
\end{displaymath}
is equal to $\frac{1000000}{1000001}$; so that there is one million to very few
to bet against one that the value $\frac{1}{1070.35}$ is not in error of
one hundredth of its value; or, which is more or less the same, that after one
century of new observations, added to the former and discussed in the same manner,
the new result does not differ from the former of more than one hundredth of its value.

Newton had found, from the observations of Pound, on the elongations of
Jupiter's satellites, the mass of this planet equal to the $1067^\textmd{th}$
part of the Sun, which differs very few from the result of Bouvard.

The mass of Uranus is $\frac{1+z}{19504}$. Replacing for $z$ its former value, 
this mass becomes $\frac{1}{17907}$. The probability that this value is between the limits
\begin{displaymath}
\frac{1}{17907} \pm \frac{1}{4}\frac{1}{19504}
\end{displaymath}
is equal to $\frac{2508}{2509}$, and the probability that this mass is between the limits
\begin{displaymath}
\frac{1}{17907} \pm \frac{1}{5}\frac{1}{19504}
\end{displaymath}
is equal to $\frac{215.6}{216.6}$.

%
%
\newpage
\indent
Les perturbations qu'Uranus produit dans le mouvement de Saturne \'etant peu
consid\'erables, on ne doit pas encore attendre des observations de ce
mouvement une grande pr\'ecision dans la valeur de sa masse. Mais, apr\`es un
si\`ecle de nouvelles observations, ajout\'ees aux pr\'ec\'edentes et
discut\'ees de la m\^eme mani\`ere, la valeur de $P$ augmentera de mani\`ere
\`a donner cette masse avec une grande probabilit\'e que sa valeur sera
contenue dans d'\'etroites limites; ce qui sera de beaucoup pr\'ef\'erable \`a
l'emploi des \'elongations des satellites d'Uranus, \`a cause de la
difficult\'e d'observer ces \'elongations.

Bouvard, en appliquant la m\'ethode pr\'ec\'edente aux $126$ \'equations de
condition que lui ont donn\'ees les observations de Jupiter et en supposant la
masse de Saturne \'egale \`a $\frac{1+z}{3534,08}$, a trouv\'e
\begin{displaymath}
z = 0,00620
\end{displaymath}
et
\begin{displaymath}
\log P = 4,8856829.
\end{displaymath}
Ces valeurs donnent la masse de Saturne \'egale \`a $\frac{1}{3512,3}$, et la
probabilit\'e que cette masse est comprise dans les limites
\begin{displaymath}
\frac{1}{3512,3} \pm \frac{1}{100}\frac{1}{3534,08}
\end{displaymath}
est \'egale \`a $\frac{11327}{11328}$.

Newton avait trouv\'e par les observations de Pound sur la plus grande
\'elongation du quatri\`eme satellite de Saturne, la masse de cette plan\`ete
\'egale \`a $\frac{1}{3012}$, ce qui surpasse d'un sixi\`eme le r\'esultat
pr\'ec\'edent. Il y a des millions de milliards \`a parier contre un que celui
de Newton est en erreur, et l'on n'en sera point surpris si l'on consid\`ere la
difficult\'e d'observer les plus grandes \'elongations des satellites de
Saturne. La facilit\'e d'observer celles des satellites de Jupiter a rendu,
comme on l'a vu, beaucoup plus exacte la valeur que Newton a conclue des
observations de Pound.

%
%
\newpage
\indent

The perturbations that Uranus induces in the movement of Saturn being negligible, 
we should not expect a great accuracy in the value of the mass from these observations of the movement.
But, after a century of new observations, added to the previous and discussed in the same manner, 
the value of $P$ increases so that the mass is given with a large probability that 
its value is contained within tight bounds; which is a lot better than using the elongations of the Uranus' satellites, 
because these elongations are difficult to observe.

Bouvard, applying the former method to $126$ equations given from the observations of Jupiter and assuming that 
the mass of Saturn is equal to $\frac{1+z}{3534.08}$, has found
\begin{displaymath}
z = 0.00620
\end{displaymath}
and
\begin{displaymath}
\log P = 4.8856829.
\end{displaymath}
These values give the mass of Saturn equal to $\frac{1}{3512.3}$, and the probability that this mass is between the limits
\begin{displaymath}
\frac{1}{3512.3} \pm \frac{1}{100}\frac{1}{3534.08}
\end{displaymath}
is equal to $\frac{11327}{11328}$.

Newton has found, from Pound's observations on the largest elongations of
Saturn's fourth satellite, that the mass of this planet is equal to
$\frac{1}{3012}$, which overestimates from than one sixth the former result.
There are millions of billions to bet against one that the one of Newton is in
error, and we should not be surprised considering the difficulty to observe 
the greatest elongations of Saturn's satellites. The easiness to observe the ones from
Jupiter's satellites has given, as we have seen, a much more exact value than the one
concluded by Newton from Pound's observations.
\newpage

\section{Comments}
\label{sec:comments}

Although Laplace presents his algorithm for two variables, in Sections~\ref{sec:oo}
and~\ref{sec:tt}, we
will assume that he is only seeking one variable and its standard deviation. This makes
explanations easier. Then in Section~\ref{sec:2var}, we generalize to two variables.

\subsection{Laplace's algorithm as a factorization algorithm}
\label{sec:oo}

The procedure used by Laplace is a variant of the Cholesky factorization of the normal equations.
We do not claim that Laplace interprets his algorithm as a factorization.
We state that, in Matrix Computation term, {\em we} can interpret Laplace's algorithm as a factorization.

In Matlab notation, if we initialize \verb M  with the lower part of $\A^T\A$ in input, his algorithm writes

\begin{verbatim}
   for k=n:-1:2,
      M(1:k-1,1:k-1) = M(1:k-1,1:k-1) - M(1:k-1,k)*M(1:k-1,k)' / M(k,k); 
   end
\end{verbatim}

The operation is a symmertic rank-1 update and Laplace only updates the lower part of the matrix \verb M  at each step.
After this operation, one obtains a {\em reverse square-root-free Cholesky factorization} of the form:
$$ \left( \A^T\A \right) =  \left( \D_M^{-1} \cdot \M\right)^T \M  ,$$
where $\D_M$ is the matrix corresponding to the diagonal of $\M$.
This is a {\em reverse} Cholesky factorization because it is a upper triangular
matrix times a lower triangular matrix as opposed to being a lower triangular
matrix times a upper triangular matrix. It is {\em square root free} because
the left factor, $\left( \D_M^{-1} \cdot \M \right)^T $, has a unit diagonal and  the right factor, $\M$, has a non-unit diagonal.
reverse Cholesky (with square roots) gives
$\D_M^{1/2}$ to each factor so that the factorization writes
$$ \left( \A^T\A \right) =  \left( \D_M^{-1/2}\cdot \M\right)^T \left( \D_M^{-1/2}\cdot \M\right) = \L^T \L .$$

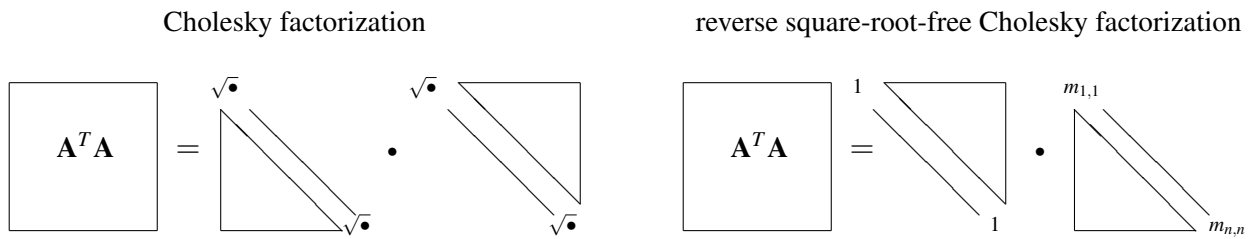
\begin{center}
\begin{tabular}{ccc}
Cholesky factorization && reverse square-root-free Cholesky factorization\\\\
\begin{picture}(220,60)(0,0)
\scriptsize

\put(  2,  2){\line(+1,+0){ 56}}
\put( 58,  2){\line(+0,+1){ 56}}
\put( 58, 58){\line(-1,+0){ 56}}
\put(  2, 58){\line(+0,-1){ 56}}

\put( 20, 30){\mbox{\normalsize$\A^T\A$}}
\put( 65, 30){\mbox{\large$=$}}

\put( 82,  2){\line(+1,+0){ 46}}
\put( 82, 48){\line(+0,-1){ 46}}
\put(128,  2){\line(-1,+1){ 46}}

\put( 78, 54){\mbox{$\sqrt{\bullet}$}}
\put(128,  2){\mbox{$\sqrt{\bullet}$}}
\put(133,  8){\line(-1,+1){ 40}}

\put(145, 30){\mbox{$\bullet$}}

\put(153, 54){\mbox{$\sqrt{\bullet}$}}
\put(206,  2){\mbox{$\sqrt{\bullet}$}}
\put(208,  8){\line(-1,+1){ 40}}

\put(172, 58){\line(+1,+0){ 46}}
\put(218, 58){\line(+0,-1){ 46}}
\put(218, 12){\line(-1,+1){ 46}}

\end{picture}
&~~~~~&
\begin{picture}(220,60)(0,0)
\scriptsize

\put(  2,  2){\line(+1,+0){ 56}}
\put( 58,  2){\line(+0,+1){ 56}}
\put( 58, 58){\line(-1,+0){ 56}}
\put(  2, 58){\line(+0,-1){ 56}}

\put( 20, 30){\mbox{\normalsize$\A^T\A$}}
\put( 65, 30){\mbox{\large$=$}}

%

\put( 66, 54){\mbox{$1$}}
\put(118,  2){\mbox{$1$}}
\put(114,  8){\line(-1,+1){ 40}}
 
\put( 78, 58){\line(+1,+0){ 46}}
\put(124, 58){\line(+0,-1){ 46}}
\put(124, 12){\line(-1,+1){ 46}}

\put(135, 30){\mbox{$\bullet$}}
  
\put(150,  2){\line(+1,+0){ 46}}
\put(150, 48){\line(+0,-1){ 46}}
\put(196,  2){\line(-1,+1){ 46}}
\put(146, 54){\mbox{$m_{1,1}$}}
\put(201,  2){\mbox{$m_{n,n}$}}
\put(201,  8){\line(-1,+1){ 40}}

%

\end{picture}
\end{tabular}
\end{center}

In Matlab notation, after Laplace's algorithm, we could
compute the solution with a backward and a forward solve. 
\begin{verbatim}
   z = (diag(diag(M))\M)' \ alpha;
   z = M \ z;
\end{verbatim}
where the first line is the backward solve and the second line is the forward solve.

In Laplace's algorithm, the backward solve is done on the fly as
\begin{verbatim}
   z = alpha;
   for k=n:-1:2,
      z(1:k-1) = z(1:k-1) - M(1:k-1,k)*z(k) / M(k,k);
   end
\end{verbatim}
On the fly means that this loop is inserted in the loop of the factorization. 
Since Laplace only wishes the first variable, the forward solve reduces to
\begin{verbatim}
   z(1) = z(1) / M(1,1);
\end{verbatim}

Another way to understand Laplace's factorization is to follow Trefethen and Bau's explanations of LU~\cite{Trefethen.Bau.97}.
Laplace is applying a sequence of unit upper triangular matrices, $\U_i$, to $\A^T\A$
in order to reduce $\A^T\A$ to lower triangular form. We obtain

$$
\left(
\begin{picture}(66,36)(2,28)
\scriptsize
\put(  6, 54){\mbox{$1$}}
\put( 12, 48){\mbox{$1$}}
\put( 18, 42){\mbox{$1$}}
\put( 46, 14){\mbox{$1$}}
\put( 52,  8){\mbox{$1$}}
\put( 58,  2){\mbox{$1$}}
\put( 44, 20){\line(-1,+1){ 20}}
\put( 58,  9){\mbox{$\star$}}
\put( 50, 30){\normalsize\mbox{$\U_{n-1}$}}
\end{picture}
\right)
\ldots
\left(
\begin{picture}(66,36)(2,28)
\scriptsize
\put(  6, 54){\mbox{$1$}}
\put( 12, 48){\mbox{$1$}}
\put( 18, 42){\mbox{$1$}}
\put( 46, 14){\mbox{$1$}}
\put( 52,  8){\mbox{$1$}}
\put( 58,  2){\mbox{$1$}}
\put( 44, 20){\line(-1,+1){ 20}}
\put( 18, 49){\mbox{$\star$}}
\put( 46, 49){\mbox{$\star$}}
\put( 52, 49){\mbox{$\star$}}
\put( 58, 49){\mbox{$\star$}}
\put( 30, 51){\line(+1,+0){ 14}}
\put( 50, 30){\normalsize\mbox{$\U_{2}$}}
\end{picture}
\right)
\left(
\begin{picture}(66,36)(2,28)
\scriptsize
\put(  6, 54){\mbox{$1$}}
\put( 12, 48){\mbox{$1$}}
\put( 18, 42){\mbox{$1$}}
\put( 46, 14){\mbox{$1$}}
\put( 52,  8){\mbox{$1$}}
\put( 58,  2){\mbox{$1$}}
\put( 44, 20){\line(-1,+1){ 20}}
\put( 12, 55){\mbox{$\star$}}
\put( 18, 55){\mbox{$\star$}}
\put( 46, 55){\mbox{$\star$}}
\put( 52, 55){\mbox{$\star$}}
\put( 58, 55){\mbox{$\star$}}
\put( 24, 57){\line(+1,+0){ 20}}
\put( 50, 30){\normalsize\mbox{$\U_{1}$}}
\end{picture}
\right)
\left(
\begin{picture}(66,36)(2,28)
\put(  5,  0){\line(+1,+0){ 60}}
\put( 65,  0){\line(+0,+1){ 60}}
\put( 65, 60){\line(-1,+0){ 60}}
\put(  5, 60){\line(+0,-1){ 60}}
\put( 25, 30){\mbox{$\A^T\A$}}
\end{picture}
\right)
=
\left(
\begin{picture}(66,36)(2,28)
\put(  5,  0){\line(+1,+0){ 60}}
\put(  5, 60){\line(+0,-1){ 60}}
\put(  5, 60){\line(+1,-1){ 60}}
\put( 15, 20){\mbox{$\M$}}
\end{picture}
\right)
$$
And with a few ``{\em strokes of luck}'', we can prove that 
$$ \left( \U_{n-1}, \ldots, \U_2, \U_1 \right)^{-1} = \left( \D_M^{-1} \cdot \M \right)^T $$
so that we finally recover our factorization:
$ \left( \A^T\A \right) =  \left( \D_M^{-1} \cdot \M\right)^T \M$.

If we follow this interpretation then, the backward solve is the application of
the sequence of unit upper triangular matrices to the right-hand side itself.

A third way to understand Laplace's factorization is to follow Laplace's
explanations.  Laplace is implicitly computing the L factor of the QL
factorization of $\A$ from the normal equations $\A^T\A$.

\subsection{Laplace's algorithm to compute the standard deviation of a variable}
\label{sec:tt}

We have seen in Section~\ref{subsubsec:zz} that we can derive quite easily the
standard deviation of the first variable of a statistical model from the QL factorization
of the regression matrix. The relation is given in Equation~(\ref{eq:mii}).

Since the L factor obtained by the QL factorization is also the reverse Cholesky
factor, Equation~(\ref{eq:mii}) explains how to compute the standard deviation of the
first variable with reverse Cholesky. We think this is the best way to
understand Laplace's algorithm with our contemporary tools.

\subsection{Laplace's algorithm to compute the standard deviation of two variables}
\label{sec:2var}

Laplace indeed is interested in the standard deviation of the two first variables. In this case,
he stops his reduction process at the 2-by-2 matrix corresponding to
$$ \A_{n-2}^T \A_{n-2} =
\left(
\begin{array}{cc}
\p_{n-2}^T \p_{n-2} & \p_{n-2}^T \q_{n-2}\\
\q_{n-2}^T \p_{n-2} & \q_{n-2}^T \q_{n-2}
\end{array}
\right)
=
\left(
\begin{array}{rr}
	 48442 &    48020 \\
	 48020 & 57725227 \\
\end{array}	  
\right).
$$
The covariance matrix of $u$ and $u'$ is then given by
\begin{eqnarray}
\nonumber \sigma_\alpha^2 \left(\A_{n-2}^T \A_{n-2}\right)^{-1} &=&
\sigma_b^2
\left(
\begin{array}{cc}
\p_{n-2}^T \p_{n-2} & \p_{n-2}^T \q_{n-2}\\
\q_{n-2}^T \p_{n-2} & \q_{n-2}^T \q_{n-2}
\end{array}
\right)^{-1}
=
\sigma_b^2
\left(
\begin{array}{rr}
	 48442 &    48020 \\
	 48020 & 57725227 \\
\end{array}	  
\right)^{-1}\\
\nonumber & = &
\sigma_b^2
\frac{1}{ 48442 * 57725227 - 48020^2 }
\left(
\begin{array}{rr}
	57725227 & - 48020 \\
	 - 48020 &   48442\\
\end{array}	  
\right).
\end{eqnarray}

If we use the fact that
 $\frac{1}{s-n}\left\|\e'\right\|^2$
is an unbiased estimate of $\sigma_b^2$ and approximate the term $s-n$ with $s$,
we get that the standard deviation of $z$, the first variable, 
and the standard deviation of $z'$, the second variable, are equal to
$$
\sigma_z^2 = 
\frac{1}{s}\left\|\e'\right\|^2
\frac{57725227}{ 48442 * 57725227 - 48020^2 }
\quad \textmd{and} \quad
\sigma_z'^2 = 
\frac{1}{s}\left\|\e'\right\|^2
\frac{48442}{ 48442 * 57725227 - 48020^2 }.
$$
Laplace gives these two formulae for the {\em poids}, $P$,  and not for the standard deviation.
Using the relation $P = ( 2 \sigma^2 )^{-1}$ we find the equation of Laplace for $z$ (top of p.\pageref{hhhhere} of this document)
\begin{displaymath}
P = \frac{s}{2S\varepsilon'^{(i)2}}\left[ 48442 - \frac{ (48020 )^2 }{57725227} \right]
\end{displaymath}
and, for $z'$,
\begin{displaymath}
P = \frac{s}{2S\varepsilon'^{(i)2}}\left[ 57725227 - \frac{ (48020 )^2 }{ 48442 } \right].
\end{displaymath}

\section{Numerical example}

There are two problem sets for Laplace to apply his algorithm. The first one
computes the mass of Jupiter, the second one computes the mass of Saturn.
Laplace uses observations from Bouvard
 where $(s=126,n=6)$ for
Jupiter and $(s=129,n=6)$ for Saturne. Actually Laplace only uses from Bouvard
 the $6$--by--$6$ normal equations and 
 the norm of the residual of the
least squares problem, $\e'$.  On the $6$ unknowns (in the $\x$ vector), Laplace only seeks one, the second variable $z'$.  The mass of
Jupiter in term of the mass of the Sun is given by $z'$ and the formula: $$
\textmd{mass of Jupiter} = \frac{1+z'}{1067.09}.$$ It turns out that the first
variable, $z$, represents the mass of Uranus through the formula $$ \textmd{mass
of Uranus} = \frac{1+z}{19504}.$$ Same approach holds for Saturn, so Laplace
will indeed compute and report the mass of three planets in his manuscript.

Note that at this time, Bouvard knew that he did not understand the behavior of
Uranus.  He conjectured that another planet should exist to explain the
anomality in the observed behavior of Uranus.  The mass of Uranus is
introduced as the auxiliary variable $z$ to try to cure the problem.  Laplace correctly predicts that
the computed mass for Uranus is not reliable. For the anecdoct, the missing planet
was Neptun and was found by Johann Gottfried Galle three years after the death
of Bouvard.

The number of operations performed by Bouvard is quite remarkable.  For the
computation of the  mass of Jupiter, Laplace accredited Bouvard for the
computation of the normal equations ($\A^T\A$) and of the residual norm
($\e' = \A \x - \b$), this makes
about $sn^2+2sn$ operations. For this numerical example, Laplace performed the
Cholesky factorization which is about $n^3/3$. This represents $6,048$
operations for Bouvard and a mere $72$ for Laplace! For the computation of the
mass of Saturn, the comparison is even worse since Bouvard performed all the
operations and reported the results to Laplace. We note that this means that Laplace
has explained his algorithm to Bouvart.

The computation of Laplace proved to be quite exact. In
Table~\ref{table:NASA}, we compare them with the current NASA values. We see
that the values for Jupiter and Saturn in Laplace are quite close from the NASA ones. The value for
Uranus is quite far as can be expected from its large variance. (Laplace would
have said its small {\em poids}.) We also note that the NASA values are within
Laplace's bounds for Saturn and Uranus.  The value for Jupiter is not within
Laplace's bound which means that the noise in the observations was not normal.

\begin{table}
\begin{center}
\begin{tabular}{|l|rcl|c|c|}
\hline
Planet 	& \multicolumn{3}{|c|}{Laplace/Bouvard} & estimation of the probability of error & NASA   \\
Jupiter&  (1059) & {\bf  1070 }&(1081)  & 1/1,000,000 & 1048  \\
Uranus &  (14564)& {\bf  17918}& (23241) & 1/2,509 & 22992  \\
Saturn &  (3477) & {\bf  3512}& (3547)  & 1/11,328 &3497  \\
\hline
\end{tabular}
\end{center}
\caption{\label{table:NASA}
Fraction of the mass of the Sun. The computed values from Bouvard are given in bold and the 
bound from Laplace in parenthesis. Laplace proved that his value for the mass of Uranus was not reliable 
(He was right.) The interval of confidence for Uranus and Saturn from Laplace are correct (i.e. the 
NASA values are in these intervals). 
}
\end{table}

In Table~\ref{table:laplace computations}, we perform the computation of
Laplace again using 64-bit arithmetic and we report the incorrect digits in his
computation. It is interesting to see that Laplace conserves a fix number of
significant digits along the computation. We can therefore say that Laplace was
computing in floating-point arithmetic.

\begin{table}
\begin{center}
\begin{tabular}{|l|}
\hline

STEP A:
\\
$
\begin{array}{rrrrrr|r}
    795938  & -12729398  &    6788.2  &  -1959.0   & 696.13  &   2602 &    7212.600 \\
            & 424865729  & -153106.5  & -39749.1   &  -5459  &   5722 & -738297.800 \\
            &            &   71.8720  &  -3.2252   & 1.2484  & 1.3371 &     237.782 \\
            &            &            &  57.1911   & 3.6213  & 1.1128 &     -40.335 \\
            &            &            &            & 21.543  & 46.310 &    -343.455 \\
            &            &            &            &         &    129 &   -1002.900 \\
\end{array}	  
$
\\
\hline
STEP B:
\\
\hline
$
\begin{array}{rrrrr|r}
    743454 &  -12844814 &    6761.23 &  -1981.45 &  -237.97 &   27441.68 \\
           &  424611920 & -153165.81 & -39798.46 & -7513.15 & -693812.58 \\
           &            &    71.8581 &   -3.2367 &   0.7684 &   248.1772 \\
           &            &            &   57.1815 &   3.2218 &   -31.6836 \\
           &            &            &           &    4.918 &    16.5783 \\
\end{array}	  
$
\\
\hline
$
\begin{array}{rrrrr|r}
    743454 &  -12844814 &    6761.23 &  -1981.45 &  -237.97 &   27441.6\color{red}4\color{black} \\
           &  424611920 & -153165.81 & -39798.46 & -7513.15 & -693812.58 \\
           &            &    71.8581 &   -3.2367 &   0.7684 &   248.1772 \\
           &            &            &   57.1815 &   3.2218 &   -31.6836 \\
           &            &            &           &    4.918 &    16.5783 \\
\end{array}	  
$
\\
\hline
STEP C:\\
\hline
$
\begin{array}{rrrr|r}
  731939.5 & -13208350 &   6798.41 & -1825.56 &   28243.85 \\
           & 413134432 & -151992.0 & -34876.7 & -668486.70 \\
           &           &   71.7381 &  -3.7401 &   245.5870 \\
           &           &           &  55.0710 &    42.5434 \\
\end{array}	  
$
\\
\hline
$
\begin{array}{rrrr|r}
  731939.\color{red}2\color{black} & -132083\color{red}60\color{black} &   6798.41 & -1825.5\color{red}5\color{black} &   28243.8\color{red}6\color{black} \\
           & 413134\color{red}201\color{black} & -15199\color{red}1.9\color{black} & -34876.\color{red}6\color{black} & -668486.\color{red}18\color{black} \\
           &        &   71.738\color{red}0\color{black} &  -3.7401 &   245.5870 \\
           &        &    &  55.07\color{red}09\color{black} &    42.54\color{red}41\color{black} \\
\end{array}	  
$
\\
\hline
STEP D:\\
\hline
$
\begin{array}{rrr|r}
  671414.7 & -14364541 &   6674.43 &  26833.55  \\
           & 391046861 & -154360.6 & -695430.0 \\
           &           &   71.4841 &  242.6977 \\
\end{array}	  
$
\\
\hline
$
\begin{array}{rrr|r}
  67142\color{red}3.6\color{black} & -143644\color{red}85\color{black} &   6674.43 &  26833.5\color{red}7\color{black}  \\
           & 39104686\color{red}9\color{black} & -154360.6 & -69542\color{red}9.6\color{black} \\
           &           &   71.4841 &  242.6977 \\
\end{array}	  
$
\\
\hline
STEP E:\\
\hline
$
\begin{array}{rr|r}
	 48442 &    48020 &   4172.95\\
	       & 57725227 & -171455.2\\
\end{array}	  
$
\\
\hline

$
\begin{array}{rr|r}
	 48\color{red}227\color{black} &    4802\color{red}1\color{black} &   417\color{red}3.00\color{black}\\
	       & 577252\color{red}58\color{black} & -171\color{red}355.9\color{black}
\end{array}	  
$
\\
\hline
STEP F:\\
\hline
$
\begin{array}{r}
	    z_0 = -0.00305\\
	    z_1 = 0.08916
\end{array}	  
$
\\
\hline
$
\begin{array}{r}
	    z_0 = -0.0030\color{red}4\color{black}\\
	    z_1 = 0.08916
\end{array}	  
$
\\
\hline
\end{tabular}
\end{center}
\caption{\label{table:laplace computations}
Laplace computations. Comparison with ``exact'' computation.
(``exact'' means 64-bit arithmetic.) First line is Laplace's value.
Second line is the value computed with 64-bit arithmetic and Laplace's algorithm.
}
\end{table}

We note that, while the condition number of $\A^T\A$ is failry large (above
$10^8$), we can equilibrate the matrix $\A^T\A$ with a diagonal scaling $\S$ equal
to the inverse of the square-root of the diagonal elements. In this case, the
scaled normal equations matrix, $\S(\A^T\A)\S$, has ones on its diagonal and its
condition number of 104.  So, up to a diagonal scaling, the system that Laplace
is considered is well-conditioned.

We can check the value given by Laplace for the variance of the variable $z$ and $z'$.
On the one hand, Laplace gives the {\em poids} of $z$ as 
	$ \log P = 2.0013595 $
so we obtain that the standard deviation of $z$ is given by
\verb 1/sqrt(2)/sqrt(10^(2.0013595))  that is
        $\sigma_z = 0.0706.$
On the other hand we can use the standard formula
  $ \sigma_z = \sigma_b \sqrt{ \textmd{entry}(1,1) \textmd{ of } \left( \A^T\A \right) ^{-1} }.$
In Matlab, this gives
\verb sqrt(31096/129)*sqrt(invATA(1,1))  we obtain
        $\sigma_z = 0.070\color{red}7.$
For the variable $z'$ Laplace gives its {\em poids} as 5.0778624, therefore
        $\sigma_z' = 0.002044343.$
Directly from the normal equations, we would have found
        $\sigma_z =  0.00204434\color{red}8.$

Laplace intreprets his result by giving an interval with a confidence level. 
For example, once, Laplace has computed $z'=0.08916$, the variable
such that the mass of Jupiter is $\frac{1+z'}{1067.09},$
and its associted {\em poids} ($P=10^{5.0778624}$), Laplace uses the fact that
$$\sqrt{\frac{P}{\pi}} \int_{-1/100}^{1/100} du~e^{-Pu^2} \approx \frac{1000000}{1000001} $$
to claim that there is one chance out of one million for the computed value of
$z'$ to be between $-1/100$ and $1/100$ of its exact value.
This means that there is one chance out of one million for the mass of
Jupiter to be between 1/1,081 and 1/1,059 the one of the Sun.

In the same manner, once, Laplace has computed $z=-0.00305$, the variable 
such that the mass of Uranus is $\frac{1+z}{17907},$
and its associted {\em poids} ($P=10^{2.0013595}$), Laplace uses the fact that
$$\sqrt{\frac{P}{\pi}} \int_{-1/4}^{1/4} du~ e^{-Pu^2} \approx \frac{2508}{2509} $$
to claim that there is one chance out of $2,509$ for the computed value
$z$ to be between $-1/4$ and $1/4$ of its exact value.
This means that there is one chance out of $2,509$ for the mass of Uranus
to be between 1/23,241 and 1/14,564 the one of the Sun.

\end{document}